\documentclass[11pt]{amsart}

\usepackage[T1]{fontenc}
\usepackage{lmodern}
\usepackage{amsmath,amssymb,amsthm,mathtools, comment}
\usepackage[margin=1.0in]{geometry}
\usepackage[colorlinks=true,linkcolor=blue,citecolor=blue,urlcolor=blue]{hyperref}

\usepackage[draft,commandnameprefix=ifneeded]{changes}
\definechangesauthor[color=red]{pe}
\definechangesauthor[color=red]{xm}
\definechangesauthor[color=red]{hx}

\numberwithin{equation}{section}

\newtheorem{prop}{Proposition}
\newtheorem{lemma}[prop]{Lemma}

\newtheorem{thm}[prop]{Theorem}
\newtheorem*{thm*}{Theorem}

\newtheorem{conj}[prop]{Conjecture}

\numberwithin{prop}{section}

\newtheorem{defn}[prop]{Definition}
\newtheorem*{defn*}{Definition}
\theoremstyle{definition}
\newtheorem{ex}[prop]{Example}
\newtheorem*{ex*}{Example}
\newtheorem{rmk}[prop]{Remark}

\theoremstyle{plain}
\newtheorem{manualtheoreminner}{Theorem}

\newenvironment{manualtheorem}[1]{%
\renewcommand\themanualtheoreminner{#1}%
\manualtheoreminner
}{\endmanualtheoreminner}

\newtheorem{manualcorollaryinner}{Corollary}
\newenvironment{manualcorollary}[1]{\renewcommand\themanualcorollaryinner{#1}\manualcorollaryinner}{\endmanualcorollaryinner}

\DeclareMathOperator{\Ric}{Ric}

\DeclareMathOperator{\Rm}{Rm}

\begin{document}

\title[]{Szeg\H{o} kernels and logarithmic flatness of K\"ahler metrics}

\begin{abstract}
Let $(L,h)$ be a positive Hermitian line bundle over a complex manifold $M$, inducing a K\"ahler metric $g$. We give a detailed account of the Hardy space and Szeg\H{o} kernel of the unit circle bundle of $(L^*,h^*)$ when $M$ is possibly noncompact. Assuming that $M$ admits a complete K\"ahler metric and that the Ricci curvature of $g$ is bounded below, we use the framework developed in our earlier work on Bergman logarithmic flatness \cite{EXX-BLF} to establish a diagonal expansion of the Szeg\H{o} kernel and characterize its logarithmic singularity in terms of the Tian--Yau--Zelditch--Catlin expansion of the Bergman kernels of $(L^k,h^k)$. This answers a question of Lu--Tian in the more general setting of possibly noncompact manifolds. We introduce the notions of Szeg\H{o} and Bergman logarithmic flatness for K\"ahler metrics, not necessarily globally polarized, and prove that every locally homogeneous K\"ahler metric on a compact complex manifold is Szeg\H{o} logarithmically flat. For K\"ahler--Einstein metrics, we prove Szeg\H{o} logarithmic flatness implies Bergman logarithmic flatness. Examples arising from Bryant's metrics show that Bergman logarithmic flatness does not imply Szeg\H{o} logarithmic flatness. We also prove that Szeg\H{o} logarithmic flatness is preserved under products and characterize when products of Bergman logarithmically flat metrics remain Bergman logarithmically flat. Combining the product formulas with known kernel computations, we show that every locally Hermitian symmetric K\"ahler metric has both properties.
\end{abstract}

\author [Ebenfelt]{Peter Ebenfelt}
\address{Department of Mathematics, University of California at San Diego, La Jolla, CA 92093, USA}
\email{{pebenfelt@ucsd.edu}}

\author[Xiao]{Ming Xiao}
\address{Department of Mathematics, University of California at San Diego, La Jolla, CA 92093, USA}
\email{{m3xiao@ucsd.edu}}

\author [Xu]{Hang Xu}
\address{School of Mathematics (Zhuhai), Sun Yat-sen University, Zhuhai, Guangdong 519082, China}
\email{{xuhang9@mail.sysu.edu.cn}}

\thanks{The third author was supported in part by the NSFC grant No. 12671105 and NSF Guangdong Province grant No. 2025A1515012071. }

\maketitle

\section{Introduction}\label{Sec Introduction}

Let $(L,h)\to M$ be a positive Hermitian line bundle over a complex manifold $M$ of dimension $n$. Let $g$ be the K\"ahler metric on $M$ corresponding to the K\"ahler form $\omega_g=\Ric(L,h)$. In this paper, we say that $(L,h)$ \emph{polarizes} the K\"ahler manifold $(M,g)$ (or simply the K\"ahler metric $g$), and call $(M,g;L,h)$ a \emph{polarized manifold}. We emphasize that $(M,g)$ is not required to be compact or complete.

Let $(L^*,h^*=h^{-1})$ be the dual bundle of $(L,h)$, and let $\pi:L^*\to M$ be the natural projection. Write
\[
D=D(L^*):=\{v\in L^*:|v|_{h^*}<1\},\qquad S=S(L^*):=\{v\in L^*:|v|_{h^*}=1\}
\]
for the disk and circle bundles of $L^*$, respectively, where $|v|_{h^*}$ denotes the norm of $v$ with respect to $h^*$. Since $(L^*,h^*)$ is negative, $S$ is strongly pseudoconvex. Throughout the paper, we use the defining function $\rho(v)=1-|v|_{h^*}^2$ for $S$.

We follow Lu--Tian \cite{LT04} to define the surface measure on $S$. Let $\alpha$ be the contact form on $S$, normalized by $\alpha(\partial_\theta)=1$, where $\partial_\theta$ is the vector field generating the $S^1$ action $r_\theta(v)=e^{i\theta}v$. We equip $S$ with the measure
\begin{equation}\label{Eq circle measure}
	d\sigma=\frac1{2\pi}\pi^*dV_g\wedge\alpha,\qquad dV_g=\frac{\omega_g^n}{n!}.
\end{equation}
In local angular coordinates, $d\sigma=(d\theta/2\pi)dV_g$. 
With our conventions, if $h=e^{-\varphi}$ in a local
holomorphic frame, then $g_{i\bar j}=\partial_i\partial_{\bar j}\varphi$ and $\omega_g=\frac{\sqrt{-1}}{2}\partial\bar\partial\varphi$.

When $M$ is compact, the Hardy space of $S$ is defined as $L^2(S)\cap\ker\overline\partial_b$, where the CR equation is understood in the sense of distributions. See more details in Zelditch \cite{Ze98}. The orthogonal projection onto this space is called the Szeg\H{o} projection. The Szeg\H{o} kernel of such circle bundles has been studied by Zelditch \cite{Ze98} and Lu--Tian \cite{LT04}, and later by many other researchers; see \cite{EnZh,ALZ13,LUZ,LMZ} and the references therein.  In this paper, we allow $M$ to be noncompact and $g$ to be incomplete. Related formulations of the Hardy space and Szeg\H{o} kernel in various noncompact settings have appeared in \cite{LMZ,LUZ}, etc. To the best of our knowledge, the analytic details of the boundary-value construction seem unavailable in the literature in this generality. For completeness, we give a detailed account of this matter in \S~\ref{Sec Szego kernel definition}.  This treatment is particularly needed for our definition of Szeg\H{o} logarithmic flatness for K\"ahler metrics, one of the main notions introduced and studied in this paper.

We briefly describe the definition here.  For $r>0$, let $S_r=\{|v|_{h^*}=r\}$ and denote by $d\sigma_r$ the pushforward of $d\sigma$ under the fiberwise dilation $v\mapsto rv$ from $S$ onto $S_r$. Define
\[
H^2(D)=\left\{f\in\mathcal O(D):\|f\|_{H^2(D)}^2:=\sup_{0<r<1}\int_{S_r}|f|^2d\sigma_r<\infty\right\}.
\]
We will show that every $f\in H^2(D)$ has a radial boundary value $f^b\in L^2(S)$ and that $f\mapsto f^b$ is a unitary map onto the closed subspace $H^2(S):=\{f^b:f\in H^2(D)\}$ of $L^2(S)$. The boundedness of point evaluations on $H^2(D)$ yields a reproducing kernel $\Pi$ such that
\[
f(x)=\int_S\Pi(x,y)f^b(y)d\sigma(y),\qquad x\in D.
\]
The kernel $\Pi$ extends to $D\times D$ as a function holomorphic in the first variable and anti-holomorphic in the second. We call $\Pi$ the Szeg\H{o} kernel associated with $H^2(S)$ and $d\sigma$. When $M$ is compact, $H^2(S)=L^2(S)\cap\ker\overline\partial_b$, and thus our definition  agrees with the usual one. 

We next recall the Bergman kernel functions of the tensor powers $(L^k,h^k)$. For $k\geq0$, denote by $A^2(M,L^k)$ the space of square-integrable holomorphic sections with respect to the inner product
\begin{equation}\label{Eq Bergman inner product}
	\langle s,s'\rangle
	=\int_M h^k(s,s')\,dV_g,
	\qquad s,s'\in A^2(M,L^k).
\end{equation}
Given any orthonormal basis $\{S_{k, i}\}_{i=1}^{N_k}$ of $A^2(M,L^k)$, the Bergman kernel function of $(L^k,h^k)$ is given by
\[
\mathcal B_k(z)=\sum_{i=1}^{N_k}|S_{k, i}(z)|_{h^k}^2.
\]
Here $N_k$ may be infinite, and $\mathcal B_k$ is understood to be $0$ if the space $A^2(M,L^k)$ is trivial.

Assume momentarily that $M$ is compact. By Boutet de Monvel--Sj\"ostrand \cite{BoSj},  the Szeg\H{o} kernel obeys the following diagonal expansion, for $v \in D,$
\begin{equation}\label{Eq compact Szego expansion}
\Pi(v,v)=\frac{\phi(v)}{\rho(v)^{n+1}}+\psi(v)\log\rho(v),\qquad \phi,\psi\in C^\infty(D\cup S).
\end{equation}
Using this Szeg\H{o} kernel expansion, Zelditch \cite{Ze98} and Catlin \cite{Ca99} independently proved the Tian--Yau--Zelditch--Catlin expansion
\begin{equation}\label{Eq compact TYZ expansion}
\mathcal B_k(z)\sim\left(\frac{k}{\pi}\right)^n\sum_{j=0}^\infty\frac{a_j(z)}{k^j}.
\end{equation}
This expansion holds in every $C^l$ norm. Lu \cite{Lu00} proved that the coefficients $a_j$ are universal polynomials in the curvature of $g$ and its covariant derivatives, and computed the first few coefficients explicitly. Many other works on Bergman kernel asymptotics can be found in \cite{Ti90,BSZ00,SZ02,Eng02,Ch03,DLM06,MM07,BBS08,HM14,Xu12,Se17,DHS24,HS22} and the references therein. This list is by no means exhaustive.

A natural question asks when the coefficient $\psi$ in \eqref{Eq compact Szego expansion} vanishes to the infinite order along $S$ (i.e., $\psi=O(\rho^\infty)$). Along this line, a pioneering work was due to Lu--Tian \cite{LT04}, where they proved that if $\psi=O(\rho^\infty)$ in \eqref{Eq compact Szego expansion}, then $a_j\equiv0$ for all $j>n$; see \cite[Theorem 2.3]{LT04}. They asked whether the converse holds \cite{LT04,Lu-PC}. In this paper, we answer this question affirmatively in a more general setting that allows $(M,g)$ to be noncompact and incomplete. While we were finalizing this paper, Yuang Shi brought to our attention the recent preprint \cite{PS26} on Lu--Tian's question for compact Sasakian manifolds. Lu--Tian \cite{LT04} also raised a question concerning the rigidity of logarithmic flatness near standard metrics. For recent developments in this direction, see \cite{LPZ26,LS26}.

We use the following version of the Bergman kernel expansion. It follows from \cite[Proposition 5.1 and Remark 5.2]{EXX-BLF} by taking the auxiliary bundle to be trivial. Write $\Ric(g)=-\frac{\sqrt{-1}}{2}\,\partial\overline\partial\log\det g$ for the Ricci form.

\begin{prop}\label{Prop TYZ expansion}
Let $(M,g;L,h)$ be a polarized manifold. Assume that both (a) and (b) hold:
\begin{enumerate}
\item[(a)] $M$ admits a complete K\"ahler metric, possibly different from $g$;
\item[(b)] $\Ric(g)\geq-C_0\omega_g$ on $M$ for some constant $C_0\geq0$.
\end{enumerate}
Then there exist unique smooth functions $a_j$ on $M$, $j\geq0$, such that
\begin{equation}\label{Eq general TYZ expansion}
\mathcal B_k(z)\sim\left(\frac{k}{\pi}\right)^n\sum_{j=0}^{\infty}\frac{a_j(z)}{k^j},\quad\text{locally uniformly for }z\in M.
\end{equation}
Moreover, each $a_j$ is given by the same universal polynomial in the curvature of $g$ and its covariant derivatives as in \eqref{Eq compact TYZ expansion}.
\end{prop}

\begin{rmk}\label{Rmk TYZ expansion meaning}
The asymptotic expansion \eqref{Eq general TYZ expansion} means that for any integers $N,l\geq0$ and any compact subset $K\Subset M$, there exists $C_{N,l,K}>0$ depending on $N,l,K$, $(M,g)$ and $(L,h)$, but independent of $k$, such that
\begin{equation}\label{Eq TYZ remainder estimates}
\left\|\mathcal B_k-\left(\frac{k}{\pi}\right)^n\sum_{j=0}^N\frac{a_j}{k^j}\right\|_{C^l(K)}\leq C_{N,l,K}k^{n-N-1},\quad \forall k\geq1.
\end{equation}
\end{rmk}

Note that (a) and (b) hold if $(M,g)$ is compact or homogeneous. The curvature assumption (b) holds if $g$ has constant Ricci eigenvalues or nonnegative Ricci curvature. In particular, it holds for K\"ahler--Einstein metrics and their products. By \cite[Theorem 1.1]{Lu00},
\[
a_0=1,\qquad a_1=\frac12R,\qquad a_2=\frac13\Delta_gR+\frac1{24}\bigl(|\Rm|^2-4|\Ric|^2+3R^2\bigr),
\]
where $\Rm$, $\Ric$ and $R$ are respectively the Riemannian, Ricci and scalar curvatures of $(M,g)$, and $\Delta_g$ is the Laplacian with respect to $g$. Since the coefficients $a_j$ are universal polynomials in the curvature of $g$ and its covariant derivatives, we can use these polynomials to define $a_j$ for an arbitrary K\"ahler metric, without any global polarization or completeness assumption.

Under the same assumptions, we prove the following expansion for the Szeg\H{o} kernel in the possibly noncompact setting.

\begin{manualtheorem}{1}\label{Thm main expansion}
Let $(M,g;L,h)$ satisfy conditions (a) and (b) in Proposition~\ref{Prop TYZ expansion}. Then there exist $\phi,\psi\in C^\infty(D\cup S)$ such that
\begin{equation}\label{Eq general Szego expansion}
\Pi(v,v)=\frac{\phi(v)}{\rho(v)^{n+1}}+\psi(v)\log\rho(v),\qquad v\in D.
\end{equation}
The functions $\phi$ and $\psi$ are unique modulo smooth functions vanishing to infinite order along $S$.
\end{manualtheorem}

We say that the Szeg\H{o} kernel $\Pi$ has no logarithmic singularity if $\psi$ in \eqref{Eq general Szego expansion}
vanishes to infinite order along $S$. For related work on the coefficients in the expansions of the Bergman and Szeg\H{o} kernels and their connections with boundary CR invariants, see \cite{Fe79,BEG94,Hi00,Hi06,Y26}. We next characterize this condition in terms of the coefficients $a_j$ in \eqref{Eq general TYZ expansion}. The key ingredient is the following theorem, which relates the coefficients $\phi$ and $\psi$ in \eqref{Eq general Szego expansion} to the coefficients $a_j$.

\begin{manualtheorem}{2}\label{Thm main coefficients}
Under the hypotheses of Theorem \ref{Thm main expansion}, there exist unique smooth functions $A_0,\ldots,A_n$ and $B_j$, $j\geq0$, on $M$ such that, near $S$,
\begin{equation}\label{Eq intro coefficient jets}
\phi=\sum_{j=0}^n A_j\rho^j+O(\rho^{n+1}),\qquad \psi=\sum_{j=0}^N B_j\rho^j+O(\rho^{N+1})
\end{equation}
for every $N\geq0$. Functions on $M$ are identified with their pullbacks under $\pi$. Each $A_j$ is an explicit linear combination of $a_0,\ldots,a_j$, and each $B_j$ is an explicit linear combination of $a_{n+1},\ldots,a_{n+1+j}$, as in \eqref{Eq explicit A coefficients} and \eqref{Eq explicit B coefficients}. Moreover, for any subset $U\subset M$ and any $m\geq0$, the following are equivalent:

(1) $B_0=\cdots=B_m=0$ on $U$;

(2) $a_{n+1}=\cdots=a_{n+1+m}=0$ on $U$.
\end{manualtheorem}

As an immediate consequence, we obtain the following corollary, which answers Lu--Tian's question \cite{LT04} in a more general setting.

\begin{manualcorollary}{1}\label{Cor SLF criterion}
Let $(M,g;L,h)$ satisfy the hypotheses of Theorem \ref{Thm main expansion}. Then the Szeg\H{o} kernel $\Pi$ has no logarithmic singularity if and only if $a_j\equiv0$ for every $j>n$. More generally, for any open subset $U\subset M$, $\psi$ vanishes to infinite order along $\pi^{-1}(U)\cap S$ if and only if $a_j|_U=0$ for every $j>n$.
\end{manualcorollary}

When $M$ is compact, this conclusion is also obtained in the recent preprint \cite[Theorem~6.6]{PS26}. The proofs in both papers follow the fiberwise summation method in \cite[\S~3]{EXX-BLF}. The preprint \cite{PS26} also studies Toeplitz coefficients and the Lu--Tian correspondence on compact Sasakian manifolds, including irregular ones. This setting includes, in particular, the circle bundles associated with compact polarized K\"ahler manifolds.

This characterization motivates the following definition for K\"ahler metrics, not necessarily globally polarized.

\begin{defn}\label{Def SLF metric}
Let $g$ be a K\"ahler metric on a complex manifold $M$ of dimension $n$, not necessarily globally polarized. We call $g$ \emph{Szeg\H{o} logarithmically flat} if, for every $p\in M$, there exist a relatively compact coordinate ball $U$ containing $p$ and a polarization $(L,h)$ of $g|_U$ such that the logarithmic coefficient of the Szeg\H{o} kernel of $S(L^*)$, defined with respect to the induced measure $d\sigma$ as in \eqref{Eq circle measure}, vanishes to infinite order along $S(L^*)$. Equivalently, $a_j\equiv0$ on $(M, g)$ for every $j>n$. Here $a_j$ denote the universal polynomials in the curvature and its covariant derivatives from Proposition~\ref{Prop TYZ expansion} (which are defined for every K\"ahler metric, without assuming that $(M,g)$ satisfies the hypotheses of Proposition~\ref{Prop TYZ expansion}.)
\end{defn}

The above definition uses only a local polarization, which always exists on a sufficiently small coordinate ball $U$: writing $\omega_g=\frac{\sqrt{-1}}2\partial\overline\partial\varphi$ for a local K\"ahler potential $\varphi$ on $U$, the trivial line bundle with Hermitian metric $h=e^{-\varphi}$ polarizes $g|_U$. Note that every coordinate ball $U$ admits a complete K\"ahler metric, and the Ricci curvature of $g$ is bounded below on $U$, provided that $U$ is sufficiently small. Corollary~\ref{Cor SLF criterion} thus justifies the equivalence in Definition~\ref{Def SLF metric}. Since the coefficients $a_j$ depend only on the local geometry of $g$, the definition is independent of the choice of $(L,h)$ and the coordinate balls. In particular, $g$ is Szeg\H{o} logarithmically flat if and only if $g|_U$ is Szeg\H{o} logarithmically flat for every open connected subset $U\subset M$. It therefore suffices to check this on an open cover of $M$.

When $M$ is compact, inspired by Lu--Tian \cite{LT04} and the corresponding result for Bergman logarithmic flatness in \cite[Corollary~1.13]{EXX-BLF},  we obtain the following consequences for Szeg\H{o} logarithmical flatness.

\begin{manualcorollary}{2}\label{Cor constant coefficients}
Let $(M,g)$ be a compact K\"ahler manifold of complex dimension $n$. Then the following are equivalent:
\begin{enumerate}
\item[(i)] $g$ is Szeg\H{o} logarithmically flat.
\item[(ii)] $a_j$ is constant on $M$ for every $j>n$.
\item[(iii)] For every $j>n$, $a_j$ is either nonnegative everywhere or nonpositive everywhere on $M$.
\end{enumerate}
\end{manualcorollary}

In the case when $(M, g)$ admits a global polarization, see also \cite[Corollary~7.3]{PS26}. Before stating the next corollary, we recall that a K\"ahler manifold $(M,g)$ is locally homogeneous if, for every pair of points $p,q\in M$, there exist neighborhoods $V_p$ and $V_q$ of $p$ and $q$, respectively, and a biholomorphism $f:V_p\to V_q$ such that $f(p)=q$ and $f^*g=g$ on $V_p$. In particular, in this case every $a_j$ is constant.

\begin{manualcorollary}{3}\label{Cor locally homogeneous}
Let $(M,g)$ be a compact K\"ahler manifold of complex dimension $n$, and assume that $(M,g)$ is locally homogeneous. Then $g$ is Szeg\H{o} logarithmically flat. Moreover, if $(L,h)$ polarizes $g$, then the unit circle bundle $S$ of the dual bundle $(L^*,h^*)$ is spherical if and only if $(M,g)$ is locally holomorphically isometric to one of the following spaces:
\begin{enumerate}
\item $(\mathbb B^n,\lambda\omega_{-1})$ for some $\lambda>0$.
\item $(\mathbb{CP}^n,\lambda\omega_1)$ for some $\lambda>0$.
\item $(\mathbb C^n,\omega_0)$.
\item $(\mathbb B^\ell\times\mathbb{CP}^{n-\ell},\lambda\omega_{-1}\times\lambda\omega_1)$ for some $1\leq\ell\leq n-1$ and $\lambda>0$.
\end{enumerate}
Here $\omega_c$ denotes the K\"ahler metric of constant holomorphic sectional curvature $c$.
\end{manualcorollary}

When $(M,g)$ is a Hermitian symmetric space of compact type, Engli\v{s}--Zhang \cite{EnZh} had already shown that the Szeg\H{o} kernel of $S$ has no logarithmic singularity, thereby providing the first examples of non-spherical compact circle bundles with this property. Arezzo--Loi--Zuddas \cite{ALZ13} extended this result to simply connected compact homogeneous polarized K\"ahler manifolds. They also proved that, in this setting, if $S$ is spherical, then $(M,g)$ is biholomorphically isometric to $(\mathbb{CP}^n,\lambda\omega_1)$ for some $\lambda>0$, as in case (2) of Corollary~\ref{Cor locally homogeneous}. Corollary~\ref{Cor locally homogeneous} further extends the vanishing result to compact locally homogeneous K\"ahler manifolds and gives a corresponding characterization of sphericity.

We next compare Szeg\H{o} logarithmic flatness with Bergman logarithmic flatness. Following \cite{EXX-BLF}, we call a germ of a strongly pseudoconvex hypersurface $\Sigma$ \emph{Bergman logarithmically flat} if there exists a smoothly bounded, strongly pseudoconvex domain $\Omega\subset\mathbb C^{n+1}$ with $\Sigma\subset\partial\Omega$ such that the coefficient function $\psi$ in the Fefferman expansion of the Bergman kernel $K_\Omega$ vanishes to infinite order along $\Sigma$. A strongly pseudoconvex hypersurface is called Bergman logarithmically flat if its germ at every point is Bergman logarithmically flat.

Recall from \cite[Proposition~1.5]{EXX-BLF} that, if $(M,g;L,h)$ is a polarized manifold of complex dimension $n$ and $M$ admits some complete K\"ahler metric, not necessarily $g$, then the Bergman kernel function $\widetilde{\mathcal B}_k$ of $L^k\otimes\mathcal C_M$ has the expansion
\begin{equation}\label{Eq twisted TYZ expansion}
\widetilde{\mathcal B}_k
\sim\frac1{\pi^n}\sum_{j=0}^\infty
\widetilde a_j k^{n-j}.
\end{equation}
Here $\mathcal C_M$ denotes the canonical line bundle of $M$, equipped with the Hermitian metric induced by $g$, and the $L^2$ inner product on sections of $L^k\otimes\mathcal C_M$ is defined using the tensor product metric and $dV_g$. Similar as $a_j$, each $\widetilde a_j$ is also a universal polynomial in the curvature of $g$ and its covariant derivatives. We use these universal polynomials to define $\widetilde a_j$ for an arbitrary K\"ahler metric, without any global polarization or completeness assumption.

We define Bergman logarithmic flatness for K\"ahler metrics as follows.
\begin{defn}\label{Def BLF metric}
Let $g$ be a K\"ahler metric on a complex manifold $M$ of dimension $n$, not necessarily globally polarized. We call $g$ \emph{Bergman logarithmically flat} if every point $p\in M$ has a relatively compact coordinate ball $U$ containing $p$ and a polarization $(L,h)$ of $g|_U$ such that the logarithmic coefficient of the Bergman kernel of $D(L^*)$ vanishes to infinite order along $S(L^*)$. Equivalently, $\widetilde a_j\equiv0$ on $M$ for every $j\geq n+2$.
\end{defn}

The equivalence in the definition follows from \cite[Theorem~1.9]{EXX-BLF}, applied on coordinate balls, each of which admits a complete K\"ahler metric. Since the invariants $\widetilde a_j$ depend only on the local geometry of $g$, the definition is independent of the choices of local polarizations and coordinate balls. In particular, $g$ is Bergman logarithmically flat if and only if $g|_U$ is Bergman logarithmically flat for every open connected subset $U\subset M$, and it suffices to check this on an open cover of $M$. If $(L,h)$ is a global polarization of $g$, then $g$ is Bergman logarithmically flat if and only if $S(L^*)$ is Bergman logarithmically flat. This follows by applying \cite[Theorem~1.12]{EXX-BLF} on coordinate balls and using the localization property of Bergman logarithmic flatness for hypersurfaces (see references in \cite{EXX-BLF}).

\begin{rmk}\label{Rmk polarization independence}
Both notions of logarithmic flatness are local properties of K\"ahler metrics. They are invariant under local holomorphic isometries, since the coefficients $a_j$ and $\widetilde a_j$ are universal polynomials in the curvature of $g$ and its covariant derivatives.
\end{rmk}

It is natural to ask how these two notions of logarithmic flatness are related.  We first prove that Szeg\H{o} logarithmic flatness implies Bergman logarithmic flatness in the K\"ahler--Einstein case.

\begin{manualtheorem}{3}\label{Thm main comparison}
	Every Szeg\H{o} logarithmically flat K\"ahler--Einstein metric is Bergman logarithmically flat.
\end{manualtheorem}

The proof of Theorem \ref{Thm main comparison} uses the K\"ahler--Einstein assumption to compare the untwisted and canonically twisted Bergman expansions in \eqref{Eq general TYZ expansion} and \eqref{Eq twisted TYZ expansion} by shifting the weight parameter. Whether Szeg\H{o} logarithmic flatness implies Bergman logarithmic flatness without this assumption remains open here.

In CR dimension one, Theorem \ref{Thm main comparison} is related to the work of Hirachi \cite{Hi93a}, Hirachi--Komatsu--Nakazawa \cite{HKN93}, and Takeuchi \cite{Ta20}.
For a relatively compact strictly pseudoconvex domain in a complex surface whose boundary admits a pseudo-Einstein contact form, Takeuchi proved that $\psi_\alpha=O(\rho^3)$ is equivalent to the chosen contact form $\alpha$ being pseudo-Einstein and the boundary being spherical. Hirachi proved the corresponding result for domains in $\mathbb C^2$ with transverse symmetry. More discussion of these results can be found in \S~\ref{Subsec KE comparison}.

We also show that Bergman logarithmic flatness does not imply Szeg\H{o} logarithmic flatness in general. 

\begin{manualtheorem}{4}\label{Thm Bryant_example}
For every $n\geq2$, there exist complete K\"ahler metrics on $\mathbb C^n$ satisfying condition (b) in Proposition~\ref{Prop TYZ expansion} that are Bergman logarithmically flat but not Szeg\H{o} logarithmically flat. Moreover, these metrics admit polarizations whose unit circle bundles are spherical.
\end{manualtheorem}

We recently learned that Lu--Shi \cite[Appendix~A]{LS26} constructed examples of compact K\"ahler manifolds whose metrics are Bergman logarithmically flat but not Szeg\H{o} logarithmically flat. We next study the behavior of logarithmic flatness under products and apply the resulting formulas to locally Hermitian symmetric K\"ahler metrics.
\begin{manualtheorem}{5}\label{Thm SLF product flatness}
	Let $(M_i,g_i)$ be K\"ahler manifolds, $i=1,2$. If both $g_1$ and $g_2$ are Szeg\H{o} logarithmically flat, then the product metric $g_1\oplus g_2$ on $M_1\times M_2$ is Szeg\H{o} logarithmically flat.
\end{manualtheorem}

\begin{manualtheorem}{6}\label{Thm BLF product flatness}
	Let $(M_i,g_i)$ be K\"ahler manifolds of complex dimensions $n_i$, $i=1,2$, and suppose that both $g_1$ and $g_2$ are Bergman logarithmically flat. Then $g_1\oplus g_2$ is Bergman logarithmically flat if and only if at least one of $\widetilde a_{n_1+1}(g_1)$ and $\widetilde a_{n_2+1}(g_2)$ vanishes identically. If both vanish identically, then $\widetilde a_j(g_1\oplus g_2)\equiv0$ for every $j>n_1+n_2$.
\end{manualtheorem}

In particular, Bergman logarithmic flatness can fail to be preserved under products. The Bryant metrics considered in \S~\ref{Subsec flatness examples} give explicit examples of this failure; see Remark~\ref{Rmk BLF product counterexample}.

We next study logarithmic flatness of locally Hermitian symmetric K\"ahler metrics. Here local Hermitian symmetry means that the curvature tensor is parallel with respect to the Levi-Civita connection. Equivalently, such a K\"ahler manifold is locally holomorphically isometric to an open subset of a simply connected (and possibly reducible) Hermitian symmetric space; see \cite{Mok89,Hel01}.
Recall that Engli\v{s}--Zhang \cite{EnZh} established the vanishing of logarithmic terms in the Szeg\H{o} kernels of circle bundles and the Bergman kernels of disk bundles associated with compact Hermitian symmetric spaces. Related kernel computations over noncompact symmetric, bounded homogeneous, and Euclidean bases were obtained by Roos \cite{Roos05}, Ahn--Park \cite{AP12}, Hao--Wang \cite{HW15}, Ishi--Park--Yamamori \cite{IPY17}, and Yamamori \cite{Yam13}. We combine these  known results with Theorems~\ref{Thm SLF product flatness} and~\ref{Thm BLF product flatness} to obtain the following corollary. It gives a unified interpretation of these earlier results in terms of logarithmic flatness of K\"ahler metrics.
\begin{manualcorollary}{4}\label{Cor Hermitian symmetric flatness}
	Every locally Hermitian symmetric K\"ahler metric is Szeg\H{o} logarithmically flat and Bergman logarithmically flat. Moreover, if $g$ is such a metric on a complex manifold of dimension $n$, then
	\begin{equation}\label{Eq symmetric coefficient vanishing}
		a_j(g)\equiv0,\qquad \widetilde a_j(g)\equiv0,\qquad j>n.
	\end{equation}
\end{manualcorollary}

The paper is organized as follows. \S~\ref{Sec Szego kernel definition} develops the Hardy space and Szeg\H{o} kernel and compares the Hardy space with the space of square-integrable CR functions. \S~\ref{Sec Szego asymptotics} proves Theorems~\ref{Thm main expansion} and~\ref{Thm main coefficients}, with the proofs of their corollaries in \S~\ref{Sec Szego corollaries}. \S~\ref{Sec BLF SLF} proves Theorems~\ref{Thm main comparison} and~\ref{Thm Bryant_example} and concludes with two conjectures. Finally, \S~\ref{Sec products symmetric} proves the product results in Theorems~\ref{Thm SLF product flatness} and~\ref{Thm BLF product flatness} and then Corollary~\ref{Cor Hermitian symmetric flatness}.

{\bf{Use of AI tools.}} The mathematical ideas and proofs were developed by the authors without the use of AI, partly building on their earlier work \cite{EXX-BLF}. During the preparation of the manuscript, the authors used OpenAI's ChatGPT to help improve the exposition and assist with the calculations in \S~\ref{Subsec flatness examples}.

\section{Definition of the Szeg\H{o} kernel}\label{Sec Szego kernel definition}

We give a detailed account of the Hardy space and Szeg\H{o} kernel for a polarized manifold $(M,g;L,h)$, without assuming compactness or completeness of $(M,g)$. We do not claim that the formulation is new. Indeed, it is a direct extension of the classical Hardy space construction on the unit disk, and related treatments have appeared in \cite{LUZ,LMZ}, etc. We provide a detailed construction in this generality for completeness and future reference. In the following subsections, we establish the existence of radial boundary values, the closedness of the resulting Hardy space, and the boundedness of point evaluations in $D$, and then construct its reproducing kernel.

\subsection{The Hardy space and its reproducing kernel}\label{Subsec Szego kernel definition}

Let $(M,g;L,h)$ be a polarized manifold of complex dimension $n$, and let $(L^*,h^*)$ be the dual Hermitian line bundle. Let $\pi:L^*\to M$ denote the natural projection. For $r>0$, set $$S_r=\{v\in L^*:|v|_{h^*}=r\},\qquad D_r=\{v\in L^*:|v|_{h^*}<r\},$$ and write $S=S_1$ and $D=D_1$. Let $\delta_r(v)=rv$ denote fiberwise dilation. In particular, $\delta_r$ maps $S$ onto $S_r$. We equip $S$ with the normalized measure $d\sigma$ as in \eqref{Eq circle measure} and $S_r$ with its pushforward $d\sigma_r=(\delta_r)_*d\sigma$. In particular,
\begin{equation}\label{Eq radial measure}
\int_{S_r}|f|^2d\sigma_r=\int_S|f\circ\delta_r|^2d\sigma.
\end{equation}
Throughout this section, $d\mu=dV_g=\omega_g^n/n!$.

\begin{defn}
The Hardy space of $D$ is
\[
H^2(D)=\left\{f\in\mathcal O(D):
\|f\|_{H^2(D)}^2:=\sup_{0<r<1}\int_{S_r}|f|^2d\sigma_r<\infty\right\}.
\]
We refer to $\|f\|_{H^2(D)}$ as the Hardy norm of $f$.
For $f\in\mathcal O(D)$ and $r\in (0, 1)$, we write $f_r:=f\circ\delta_r|_S$.
\end{defn}

Fix a trivializing coordinate chart $(U,z)$ and dual holomorphic frames $e_L$ and $e_{L^*}$. In this trivialization, write $h(z)=h(e_L,e_L)(z)$. A vector in $L^*|_U$ has the form $\xi e_{L^*}(z)$, and $$D_r\cap\pi^{-1}(U)=\{(z,\xi):|\xi|^2<r^2 h(z)\},\qquad S_r\cap\pi^{-1}(U)=\{(z,\xi):|\xi|^2=r^2h(z)\}.$$ 
In these coordinates, the dilation $\delta_r$ is given by $\delta_r(z,\xi)=(z,r\xi)$.
We use the standard local angular coordinates $(z,\theta)\mapsto(z,rh(z)^{1/2}e^{i\theta})$ on $S_r$, in which
\[
d\sigma_r=\frac{d\theta}{2\pi}d\mu,\qquad
f_r(z,\theta)=f(z,rh(z)^{1/2}e^{i\theta}).
\]
Set
\[
A_{r,f}(z)=\frac1{2\pi}\int_0^{2\pi}|f_r(z,\theta)|^2d\theta.
\]
Since $A_{r,f}(z)$ is the circular mean of $|f|^2$ over the fiber $S_r\cap L_z^*$, the function $A_{r,f}$ is independent of the choice of local coordinates and frame, and hence is globally defined on $M$.
Moreover,
\begin{equation}\label{Eq circular mean integral}
\int_{S_r}|f|^2d\sigma_r=\int_M A_{r,f}\,d\mu.
\end{equation}
For later use, note that if $f\in H^2(D)$, then
\begin{equation}\label{Eq local radial bound}
\int_U A_{r,f}\,d\mu
=\int_{S_r\cap\pi^{-1}(U)}|f|^2d\sigma_r
\leq\|f\|_{H^2(D)}^2.
\end{equation}

Let us recall the following monotonicity property for circular means of subharmonic functions. We include a proof for the reader's convenience.
\begin{lemma}\label{Lemma circular means}
Let $u$ be a continuous subharmonic function on $\{\xi\in\mathbb C:|\xi|<R\}$. Then, for $0\leq s\leq r<R$,
\[
\frac1{2\pi}\int_0^{2\pi}u(se^{i\theta})d\theta
\leq\frac1{2\pi}\int_0^{2\pi}u(re^{i\theta})d\theta.
\]
\end{lemma}

\begin{proof}
We may assume $s<r$. Let $v$ be the harmonic function on $\{|\xi|<r\}$ with boundary values $u$ on $\{|\xi|=r\}$. The maximum principle gives $u\leq v$ on $\{|\xi|<r\}$. Hence the mean-value property implies
\[
\frac1{2\pi}\int_0^{2\pi}u(se^{i\theta})d\theta
\leq\frac1{2\pi}\int_0^{2\pi}v(se^{i\theta})d\theta
=v(0)=\frac1{2\pi}\int_0^{2\pi}u(re^{i\theta})d\theta.
\]
\end{proof}

Applying Lemma~\ref{Lemma circular means} to $|f(z,\cdot)|^2$, we see that $A_{r,f}(z)$ is nondecreasing in $r$. Integrating over $M$ and using \eqref{Eq circular mean integral} together with \eqref{Eq radial measure}, we obtain
\begin{equation}\label{Eq Hardy norm limit}
\|f\|_{H^2(D)}^2
=\lim_{r\to1^-}\int_{S_r}|f|^2d\sigma_r
=\lim_{r\to1^-}\|f_r\|_{L^2(S)}^2.
\end{equation}

\begin{prop}\label{Prop radial boundary values}
Every $f\in H^2(D)$ has a radial boundary value $f^b\in L^2(S)$. More precisely,
\[
f(rv)\longrightarrow f^b(v)\quad\text{for almost every }v\in S,
\qquad
f_r\longrightarrow f^b\quad\text{in }L^2(S)
\]
as $r\to1^-$. Moreover,
\begin{equation}\label{Eq boundary isometry}
\|f^b\|_{L^2(S)}=\|f\|_{H^2(D)}.
\end{equation}
\end{prop}

\begin{proof}
Let $f\in H^2(D)$. Since $A_{r,f}$ is nonnegative and nondecreasing in $r$, the monotone convergence theorem, together with \eqref{Eq circular mean integral} and \eqref{Eq Hardy norm limit}, yields
\[
\int_M\lim_{r\to1^-}A_{r,f}(z)\,d\mu
=\lim_{r\to1^-}\int_M A_{r,f}(z)\,d\mu
=\|f\|_{H^2(D)}^2<\infty.
\]
Consequently, $\lim_{r\to1^-}A_{r,f}(z)<\infty$ for almost every $z\in M$. For each such $z$, the restriction of $f$ to the fiber disk belongs to the classical Hardy space. By the classical Fatou theorem, it therefore has radial limits almost everywhere on the fiber circle.

Let $$\mathcal T=\{v\in S:\lim_{r\to1^-}f(rv)\text{ exists in }\mathbb C\}.$$ 
We denote by $\sigma$ the Borel measure on $S$ induced by the volume form $d\sigma$ in \eqref{Eq circle measure}. By using continuity of $f(rv)$ in $r \in (0, 1),$ one can readily show $\mathcal T$ is a $\sigma-$measurable subset of $S$.
Fix a local trivialization of $L$ and $L^*$ over an open set $U$. Tonelli's theorem gives
\[
\sigma\big((S\setminus\mathcal T)\cap\pi^{-1}(U)\big)
=\int_U\left(\frac1{2\pi}\int_0^{2\pi}\mathbf1_{S\setminus\mathcal T}(z,\theta)\,d\theta\right)d\mu(z)=0,
\]
since the inner integral is zero for almost every $z\in U$. Applying this argument to all $U$ in a countable trivializing cover of $M$, we obtain $\sigma(S\setminus\mathcal T)=0$.
Thus, we can define the radial limit:
\[
f^b(v):=\lim_{r\to1^-}f(rv),\qquad\text{for almost every }v\in S.
\]
The function $f^b$ is measurable. By Fatou's lemma and \eqref{Eq Hardy norm limit}, $\|f^b\|_{L^2(S)}\leq\|f\|_{H^2(D)}$. Hence $f^b \in L^2(S)$.

To prove the $L^2$ convergence assertion, we fix a local trivialization over $U$ as above and write the fiberwise Taylor expansion of $f$ as
\begin{equation}\label{Eq section2 fiber expansion}
f(z,\xi)=\sum_{m=0}^\infty b_m(z)\xi^m.
\end{equation}
Here each $b_m$ is holomorphic, and the series converges uniformly on compact subsets of $D\cap\pi^{-1}(U)$. Note that $|b_m|^2h^m$ is independent of the choice of frame. We may thus set $\lambda_m=\int_M|b_m|^2h^m dV_g$.  Applying Parseval's identity on each fiber, we obtain
\[
\frac1{2\pi}\int_0^{2\pi}|f_r(z,\theta)|^2d\theta=\sum_{m=0}^\infty|b_m(z)|^2h(z)^mr^{2m}.
\]
Integrating over $M$ and applying Tonelli's theorem yields the first equality below. The second follows by letting $r\to1^-$ and using monotone convergence and \eqref{Eq Hardy norm limit}:
\begin{equation}\label{Eq global coefficient norm}
\|f_r\|_{L^2(S)}^2=\sum_{m=0}^\infty\lambda_mr^{2m},\qquad \sum_{m=0}^\infty\lambda_m=\|f\|_{H^2(D)}^2<\infty.
\end{equation}
For $r,t<1$, we also have $\|f_r-f_t\|_{L^2(S)}^2=\sum_m\lambda_m|r^m-t^m|^2$. Dominated convergence shows that this tends to zero as $r,t\to1^-$. Thus $f_r$ converges in $L^2(S)$ as $r\to1^-$. Since $L^2$ convergence implies almost-everywhere convergence along a subsequence, and $f_r$ converges to $f^b$ almost everywhere on $S$, the $L^2$ limit must coincide with $f^b$ almost everywhere. Hence $f_r\to f^b$ in $L^2(S)$. Combining this with \eqref{Eq Hardy norm limit}, we obtain \eqref{Eq boundary isometry}.
\end{proof}

We are now ready to define the Hardy space of the circle bundle $S$.

\begin{defn}\label{Def Hardy space on S}
The Hardy space on $S$ is $H^2(S)=\{f^b:f\in H^2(D)\}\subset L^2(S)$.
\end{defn}

Proposition~\ref{Prop radial boundary values} shows that the boundary map $b:f\mapsto f^b$ is a linear isometry onto $H^2(S)$. In particular, $b$ is injective. We equip $H^2(D)$ with the inner product induced from $H^2(S)$:
$$\langle f,g\rangle_{H^2(D)}:=\langle f^b,g^b\rangle_{L^2(S)}.$$
We next prove the completeness of $H^2(S)$ and $H^2(D)$. To this end, we first establish the following estimate.
\begin{lemma}\label{Lemma Hardy evaluation}
For every compact subset $K\Subset D$, there is a constant $C_K>0$ such that
\[
\sup_{v\in K}|f(v)|\leq C_K\|f\|_{H^2(D)}=C_K\|f^b\|_{L^2(S)},\qquad f\in H^2(D).
\]
\end{lemma}

\begin{proof}
Choose $0<r_0<1$ such that $K\subset D_{r_0}$. Fix $p\in K$, and choose a trivialization $(z,\xi)$ near $p=(z_0,\xi_0)$. Within the coordinate chart at $z_0$, choose $V\Subset W\Subset U$, with $z_0\in V$, and a disk $\Delta_\varepsilon=\{\xi:|\xi-\xi_0|<\varepsilon\}$ such that $W\times\Delta_\varepsilon\subset D_{r_0}\cap\pi^{-1}(U)$. For $r_0<r<1$, let $\Gamma_r(z)=\{\eta:|\eta|=rh(z)^{1/2}\}$. Cauchy's integral formula gives
\[
f(z,\xi)=\frac1{2\pi i}\int_{\Gamma_r(z)}\frac{f(z,\eta)}{\eta-\xi}d\eta,
\qquad (z,\xi)\in W\times\Delta_\varepsilon.
\]
Writing $\eta=rh(z)^{1/2}e^{i\theta}$ and using $|\xi|<r_0h(z)^{1/2}$, we have $|\eta-\xi|\geq(r-r_0)h(z)^{1/2}$. Using this and the Cauchy--Schwarz inequality yields, for $(z,\xi)\in W\times\Delta_\varepsilon,$
\begin{align*}
|f(z,\xi)|
\leq\frac{r}{r-r_0}\int_0^{2\pi}|f_r(z,\theta)|\frac{d\theta}{2\pi}
\leq\frac1{r-r_0}
\left(\int_0^{2\pi}|f_r(z,\theta)|^2\frac{d\theta}{2\pi}\right)^{1/2}.
\end{align*}
Squaring the above bound for $|f(z,\xi)|$ and integrating over $W$, we obtain, by \eqref{Eq local radial bound},
\[
\int_W|f(z,\xi)|^2d\mu
\leq\frac1{(r-r_0)^2}\|f\|_{H^2(D)}^2.
\]
Letting $r\to1^-$ gives
\[
\int_W|f(z,\xi)|^2d\mu
\leq\frac1{(1-r_0)^2}\|f\|_{H^2(D)}^2,
\qquad \xi\in\Delta_\varepsilon.
\]
By the sub-mean value inequality for holomorphic functions and the comparability of $d\mu$ with Euclidean volume on $\overline W$, there is a constant $C_{V,W}$ such that $\sup_{z\in V}|f(z,\xi)| \leq C_{V,W}\|f(\cdot,\xi)\|_{L^2(W,d\mu)}$. Combining these estimates, we obtain
\[
\sup_{(z,\xi)\in V\times\Delta_\varepsilon}|f(z,\xi)|
\leq \frac{C_{V,W}}{1-r_0}\|f\|_{H^2(D)}.
\]
By compactness, we can cover $K$ by finitely many such neighborhoods. The proof is therefore completed by taking the maximum of the corresponding constants.
\end{proof}

\begin{prop}\label{Prop Hardy closedness}
The space $H^2(S)$ is closed in $L^2(S)$. Consequently, $H^2(D)$ and $H^2(S)$ are both Hilbert spaces, and the boundary map $b:H^2(D)\to H^2(S)$ is unitary.
\end{prop}

\begin{proof}
Fix any Cauchy sequence $\{f_i^b\}$ in $H^2(S)$. Let $f_i\in H^2(D)$ be the unique functions with boundary values $f_i^b$. Then, by \eqref{Eq boundary isometry}, $\{f_i\}$ is Cauchy in the Hardy norm. Applying Lemma~\ref{Lemma Hardy evaluation} to $f_i-f_k$, we see that $\{f_i\}$ is uniformly Cauchy on every compact subset of $D$. Its limit $f_0$ is therefore holomorphic on $D$.
For fixed $0<r<1$, we have $f_i(rv)\to f_0(rv)$ for every $v\in S$. Fatou's lemma gives
\[
\int_S|f_0(rv)|^2d\sigma(v)
\leq\liminf_{i\to\infty}\int_S|f_i(rv)|^2d\sigma(v)
\leq\sup_i\|f_i\|_{H^2(D)}^2<\infty.
\]
Since the bound is independent of $r$, we have $f_0\in H^2(D)$. We next show that $f_i\to f_0$ in the Hardy norm. Fix $\varepsilon>0$ and choose $N$ such that $\|f_i-f_k\|_{H^2(D)}<\varepsilon$ for any $i,k\geq N$. Applying Fatou's lemma once more, for any $i\geq N$ and $0<r<1$ we obtain
\[
\int_S|f_i(rv)-f_0(rv)|^2d\sigma(v)
\leq\liminf_{k\to\infty}\int_S|f_i(rv)-f_k(rv)|^2d\sigma(v)
\leq\varepsilon^2.
\]
Taking the supremum in $r$ gives $\|f_i-f_0\|_{H^2(D)}\leq\varepsilon$ for any $i\geq N$. This shows that $f_i\to f_0$ in $H^2(D)$. By \eqref{Eq boundary isometry}, we also have $f_i^b\to f_0^b$ in $L^2(S)$. Therefore, $H^2(S)$ is closed in $L^2(S)$ and hence is a Hilbert space. By the discussion after Definition~\ref{Def Hardy space on S}, $H^2(D)$ is also a Hilbert space and the boundary map $b:H^2(D)\to H^2(S)$ is unitary.
\end{proof}

We now define the Szeg\H{o} kernel. For $x\in D$, the evaluation map $$P_x:H^2(S)\longrightarrow\mathbb C,\qquad P_x(f^b)=f(x),$$ is well-defined and bounded by Lemma~\ref{Lemma Hardy evaluation}. By the Riesz representation theorem, there exists a unique $g_x\in H^2(S)$ such that $f(x)=\langle f^b,g_x\rangle_{L^2(S)}$. We define $\Pi(x,\cdot)=\overline{g_x}\in L^2(S)$. Thus
\begin{equation}\label{Eq Szego boundary reproduction}
f(x)=\int_S\Pi(x,y)f^b(y)d\sigma(y),\qquad f\in H^2(D).
\end{equation}

\begin{prop}\label{Prop Szego basis expansion}
Let $\{\psi_i\}_{i=1}^N$, $0\leq N\leq\infty$, be an orthonormal basis of $H^2(S)$, and let $\Psi_i\in H^2(D)$ satisfy $\Psi_i^b=\psi_i$. Then, for each $x\in D$,
\begin{equation}\label{Eq Szego mixed expansion}
\Pi(x,\cdot)=\sum_{i=1}^N\Psi_i(x)\overline{\psi_i}
\quad\text{in }L^2(S).
\end{equation}
The formula
\begin{equation}\label{Eq Szego interior expansion}
\Pi(x,y)=\sum_{i=1}^N\Psi_i(x)\overline{\Psi_i(y)},\qquad x,y\in D,
\end{equation}
defines a function holomorphic in $x$ and anti-holomorphic in $y$. The series converges uniformly on compact subsets of $D\times D$. For each fixed $y\in D$, $\Pi(\cdot,y)$ belongs to $H^2(D)$ and has radial boundary value
\begin{equation}\label{Eq Szego radial boundary expansion}
\Pi(\cdot,y)^b=\sum_{i=1}^N\psi_i\overline{\Psi_i(y)}
\quad\text{in }H^2(S).
\end{equation}
Similarly, for fixed $x\in D$, $\overline{\Pi(x,\cdot)}$ belongs to $H^2(D)$ and its boundary value is $g_x$.
\end{prop}

\begin{proof}
For fixed $x\in D$, the defining property of $g_x$ gives $\Psi_i(x)=\langle\psi_i,g_x\rangle_{L^2(S)}$. We expand $g_x$ in the orthonormal basis to obtain
\[
g_x=\sum_i\langle g_x,\psi_i\rangle_{L^2(S)}\psi_i=\sum_i\overline{\Psi_i(x)}\psi_i,
\]
with convergence in $H^2(S)$. This proves \eqref{Eq Szego mixed expansion}. For $x\in K\Subset D$, let $f\in H^2(D)$ have boundary value $g_x$. Then $f(x)=\langle f^b,g_x\rangle_{L^2(S)}
=\|g_x\|_{L^2(S)}^2$ and $\|f\|_{H^2(D)}=\|g_x\|_{L^2(S)}$. Applying Lemma~\ref{Lemma Hardy evaluation} to $f$ at $x$, we obtain
\[
f(x)=\|g_x\|_{L^2(S)}^2
\leq C_K\|g_x\|_{L^2(S)},
\]
which implies $\|g_x\|_{L^2(S)}\leq C_K$. Therefore,
\begin{equation}\label{Eq Szego diagonal bound}
\sum_i|\Psi_i(x)|^2=\|g_x\|_{L^2(S)}^2
\leq C_K^2,\qquad x\in K.
\end{equation}
In particular, for each fixed $y\in D$, we have $\sum_i|\Psi_i(y)|^2<\infty$. Since $\{\Psi_i\}$ is an orthonormal basis of $H^2(D)$, the series $\sum_i\Psi_i(\cdot)\overline{\Psi_i(y)}$ converges in $H^2(D)$. By Lemma~\ref{Lemma Hardy evaluation}, the convergence in $H^2(D)$ implies locally uniform convergence on $D$. Hence the sum $\sum_i\Psi_i(x)\overline{\Psi_i(y)}$ is holomorphic in $x$. Similarly, for each fixed $x\in D$, the sum is anti-holomorphic in $y$. By Hartogs' theorem, the sum is holomorphic in $(x,\overline y)$.
In particular, $\Pi(x,x)=\sum_{i}|\Psi_i(x)|^2$ is continuous. By Dini's theorem, the convergence is uniform on compact subsets of $D$.
Then Cauchy--Schwarz inequality shows local uniform convergence of \eqref{Eq Szego interior expansion}. Finally, the boundary identity \eqref{Eq Szego radial boundary expansion} follows from the continuity of the unitary map $f\mapsto f^b$. The last statements about $\overline{\Pi(x,\cdot)}$ can be proved similarly.
\end{proof}

The function $\Pi$ is independent of the orthonormal basis and satisfies $\Pi(x,y)=\overline{\Pi(y,x)}$ for all $x,y\in D$. We call $\Pi$ the Szeg\H{o} kernel associated with $H^2(S)$ and the measure $d\sigma$.

\subsection{Comparison with the space of square-integrable CR functions}\label{Subsec Szego kernel comparison}
In the classical setting, the base manifold $M$ is assumed to be compact, and the Hardy space of the circle bundle is usually defined by $\widehat H^2(S)=L^2(S)\cap\ker\overline\partial_b$, where the CR equation is understood in the distribution sense. See more details in \cite{Ze98} and \cite{LT04}, where the space is also described in terms of boundary values of holomorphic functions on $D$. Although this definition still makes sense when $M$ is noncompact, functions in $\widehat H^2(S)$ need not be boundary values of holomorphic functions on $D$, as Example~\ref{E1} below shows. We therefore define the Hardy space using radial boundary values and will discuss its relation to $\widehat H^2(S)$ in this subsection over a general (not necessarily compact) base manifold $M$. We first prove the following result.

\begin{prop}\label{Prop comparison Hardy spaces}
It always holds that $H^2(S)\subseteq\widehat H^2(S)$. If $M$ is compact, then $H^2(S)=\widehat H^2(S)$.
\end{prop}

The assertion in the compact case is standard. We include a proof here for the reader's convenience. Let $\overline\partial_b$ denote the tangential Cauchy--Riemann operator on $S$. Set $H^\infty(S)=\{u\in C^\infty(S):\overline\partial_bu=0\}$. We first recall the following standard consequence of Fourier decomposition.

\begin{lemma}\label{Lemma smooth CR density}
If $M$ is compact, then $H^\infty(S)$ is dense in $\widehat H^2(S)$ in the $L^2$ norm.
\end{lemma}

If $M$ is compact, so is $S$. In this case, the Kohn--Rossi extension theorem \cite{KR65} implies that every $u\in C^\infty(S)$ that is CR extends to a function $F\in\mathcal O(D)\cap C^\infty(\overline D)$.

\begin{proof}[Proof of Proposition~\ref{Prop comparison Hardy spaces}]
First suppose $M$ is compact. For $u\in H^\infty(S)$, let $F\in\mathcal O(D)\cap C^\infty(\overline D)$ be the extension just described. By compactness, $F_r\to u$ uniformly on $S$ and $\sup_{0<r<1}\|F_r\|_{L^2(S)}<\infty$. Thus $F\in H^2(D)$ and $F^b=u$, so $H^\infty(S)\subseteq H^2(S)$. Lemma~\ref{Lemma smooth CR density} and Proposition~\ref{Prop Hardy closedness} imply $\widehat H^2(S)=\overline{H^\infty(S)}^{\,L^2} \subseteq H^2(S)$.

We now prove $H^2(S)\subseteq\widehat H^2(S)$ without assuming compactness. If $f^b\in H^2(S)$, choose $f\in H^2(D)$ with radial boundary value $f^b$. For every $0<r<1$, $f_r=f\circ\delta_r|_S$ is smooth and CR. Since $f_r\to f^b$ in $L^2(S)$ and $\overline\partial_b f_r=0$, we have $\overline\partial_b f^b=0$ in the distribution sense.  This proves $H^2(S)\subseteq\widehat H^2(S)$.
\end{proof}

\begin{rmk}\label{Rmk complete Hardy spaces}
We thank Matteo Palmieri and Wei-Chuan Shen for pointing out that $H^2(S)=\widehat H^2(S)$ when $g$ is complete. In this case, the orthogonal projection onto $H^2(S)$ agrees with the Szeg\H{o} projection considered in \cite{HMW22}. If, in addition, the Ricci curvature of $g$ is bounded below relative to $g$, their Theorem~1.2 shows that this projection is a Fourier integral operator with complex phase. Under the same assumptions, their Theorem~1.6 gives a full asymptotic expansion of the Bergman kernels of $L^k$ as $k\to\infty$.
\end{rmk}

We conclude this section with an example of a noncompact manifold $M$ with an incomplete metric $g$ for which $H^2(S)\neq\widehat H^2(S)$, in contrast to the cases discussed in Proposition~\ref{Prop comparison Hardy spaces} and Remark~\ref{Rmk complete Hardy spaces}.

\begin{ex}\label{E1}
Let $M$ be the unit disk $\mathbb{D}=\{z\in\mathbb C:|z|<1\}$, and $L=M\times\mathbb C$ be the trivial line bundle equipped with the Hermitian metric $h=e^{-|z|^2}$. Then $(L,h)$ induces the Euclidean metric on $M$. In this setting,
\[
D=\{(z,\xi):|z|<1,\ |\xi|^2e^{|z|^2}<1\},\qquad S=\{(z,\xi):|z|<1,\ |\xi|^2e^{|z|^2}=1\}.
\]
The function $F=\xi^{-1}|_S$ is smooth and CR. In angular coordinates, $F(z,\theta)=e^{|z|^2/2}e^{-i\theta}$, so
\[
\|F\|_{L^2(S)}^2=\int_{|z|<1}e^{|z|^2}dV_g<\infty.
\]
Thus $F\in\widehat H^2(S)$.
For every $f\in H^2(D)$, holomorphicity in the fiber variable gives $\int_0^{2\pi}f_r(z,\theta)e^{i\theta}\,d\theta=0$ for every $z\in M$ and $0<r<1$. By the Cauchy--Schwarz inequality,
\[
\begin{aligned}
\int_M\left|\frac1{2\pi}\int_0^{2\pi}f^b(z,\theta)e^{i\theta}\,d\theta\right|^2dV_g
&\leq\int_M\frac1{2\pi}\int_0^{2\pi}|f^b(z,\theta)-f_r(z,\theta)|^2\,d\theta\,dV_g\\
&=\|f^b-f_r\|_{L^2(S)}^2\to0.
\end{aligned}
\]
Thus the angular integral of $f^b$ vanishes for almost every $z$. However,
\[
\frac1{2\pi}\int_0^{2\pi}F(z,\theta)e^{i\theta}\,d\theta=e^{|z|^2/2}\ne0.
\]
Hence $F\notin H^2(S)$. This shows $H^2(S)\subsetneq\widehat H^2(S)$
\end{ex}

\section{The Szeg\H{o} kernel and its asymptotic expansion}\label{Sec Szego asymptotics}

\subsection{The relation between the Szeg\H{o} kernel and the Bergman kernels}\label{Subsec Szego Bergman relation}

Let $(M,g;L,h)$ be a polarized manifold, and let $(L^*,h^*)$ be the dual line bundle of $(L,h)$. Write $\pi:L^*\to M$ for the natural projection and $D=D(L^*)=\{v\in L^*: \rho(v):=1-|v|_{h^*}^2>0\},\, S=S(L^*)=\partial D$. We use the measure $d\sigma$ and the Hardy space $H^2(S)$ defined in \S~\ref{Sec Szego kernel definition}. For every integer $m\geq0$, write $\mathcal B_m$ for the Bergman kernel function of $(L^m,h^m)$, with respect to the volume form $dV_g$. Here $L^0$ is the trivial line bundle with its constant metric. The following relation is more or less well-known to experts in many settings (see \cite{Ze98}; see also \cite{LUZ, LMZ} for a related formulation in the noncompact setting). With our definitions of the Hardy space and Szeg\H{o} kernel, the same argument applies without assuming compactness of $M$. We include the proof following \cite{Ze98} for completeness and convenience of readers.

\begin{thm}\label{Thm Szego Bergman relation}
The Szeg\H{o} kernel $\Pi$ of $S$ and the Bergman kernel functions $\mathcal B_m$ satisfy
\begin{equation}\label{Eq Szego Bergman series}
\Pi(v,v)=\sum_{m=0}^{\infty}\mathcal B_m(\pi(v))(1-\rho(v))^m,\qquad v\in D.
\end{equation}
The infinite series above converges uniformly on every compact subset of $D$.
\end{thm}

In \eqref{Eq Szego Bergman series}, $(1-\rho)^0$ is understood to be $1$, including on the zero section $M\subset L^*$.
Following \cite{Ze98}, we decompose the Hardy space into weight spaces under the $S^1$-action on $S$. Choose a local coordinate chart $(U,z)$ of $M$ and a holomorphic frame $e_{L^*}$ of $L^*$ over $U$. This induces a trivialization $\pi^{-1}(U)\ni \xi e_{L^*}(z)\longmapsto(z,\xi)\in U\times\mathbb C$. Let $e_L$ be the dual frame of $e_{L^*}$,
and write $h(z)=h(e_L(z),e_L(z))$ for simplicity.
Consider the global $S^1$ action $r_\theta(v)=e^{i\theta}v$ on $L^*$, which in the above trivialization is given by $r_\theta(z,\xi)=(z,e^{i\theta}\xi)$. Since $r_\theta$ preserves $D$ and $d\sigma$, the pullback map $r_\theta^*$ defines a unitary action on $H^2(S)$. For $m\in\mathbb Z$, set 
\[H_m^2(S)=\{u\in H^2(S):r_\theta^*u=e^{im\theta}u\text{ for every }\theta\in\mathbb R\}.\]
Then each $H_m^2(S)$ is a closed subspace of $H^2(S)$.

\begin{prop}\label{Prop Hardy Fourier decomposition}
We have $H_m^2(S)=\{0\}$ for $m<0$, and the orthogonal decomposition
\begin{equation}\label{Eq Hardy Fourier decomposition}
H^2(S)=\bigoplus_{m=0}^{\infty}{}^\perp H_m^2(S).
\end{equation}
\end{prop}

\begin{proof}
For $m\in\mathbb Z$, the Fourier projection on $L^2(S)$ is given by
\begin{equation}\label{Eq Fourier projection}
	(P_m u)(v)=\frac1{2\pi}\int_0^{2\pi}u(r_\theta v)e^{-im\theta}\,d\theta.
\end{equation}
Let $f\in H^2(D)$ and write $f_r(v)=f(rv)$ for $v\in S$ and $0<r<1$. Locally, holomorphicity in the fiber variable gives
\begin{equation}\label{Eq Hardy fiber Taylor}
f(z,\xi)=\sum_{j=0}^{\infty}b_j(z)\xi^j,
\end{equation}
where the coefficient functions $b_j$ are holomorphic on $U$. Thus $P_m f_r=0$ for $m<0$. Since $f_r\to f^b$ in $L^2(S)$ as $r\to 1^-$ by Proposition \ref{Prop radial boundary values}, we also have $P_m f^b=0$ for $m<0$.

	For $m\geq0$, define the holomorphic function
\[
F_m(v)=\frac{1}{2\pi}\int_0^{2\pi}
f(r_\theta v)e^{-im\theta}\,d\theta,
\qquad v\in D.
\]
Then $(F_m)_r=P_m f_r$. Since $P_m$ is an orthogonal
projection on $L^2(S)$, we have
\[
\|(F_m)_r\|_{L^2(S)}
=\|P_m f_r\|_{L^2(S)}
\leq\|f_r\|_{L^2(S)}.
\]
Taking the supremum over $0<r<1$, we obtain
\[
\|F_m\|_{H^2(D)}\leq\|f\|_{H^2(D)}.
\]
In particular, $F_m\in H^2(D)$. By Proposition~2.3
and the boundedness of $P_m$, we have
\[
(F_m)_r=P_m f_r\longrightarrow P_m f^b
\quad\text{in }L^2(S)\quad\text{as }r\to1^-.
\]
Moreover, a change of variables shows that $r_\theta^*F_m=e^{im\theta}F_m$. Since the boundary map commutes with the $S^1$ action, we obtain
\[
P_m f^b=(F_m)^b\in H_m^2(S).
\]
Conversely, if $u\in H_m^2(S)$, then $P_m u=u$ by \eqref{Eq Fourier projection}. Thus $P_m$ maps $H^2(S)$ onto
$H_m^2(S)$. In particular, $H_m^2(S)=\{0\}$ for $m<0$.

Finally, for each $z\in M$, let
$S_z=\pi^{-1}(z)\cap S$ denote the circle fiber over $z$.
By Fubini's theorem, for almost every $z\in M$,
the restriction of $f^b$ to $S_z$ belongs to $L^2(S_z)$.
Since $P_m f^b=0$ for $m<0$, the classical Fourier series
theorem and Parseval's identity give
\[
f^b=\sum_{m=0}^{\infty}P_m f^b
\quad\text{in }L^2(S_z),\qquad
\|f^b\|_{L^2(S_z)}^2
=\sum_{m=0}^{\infty}\|P_m f^b\|_{L^2(S_z)}^2
\]
for almost every $z\in M$. By the orthogonality of
the Fourier components,
\[
\left\|f^b-\sum_{m=0}^{N}P_m f^b
\right\|_{L^2(S_z)}^2
\leq \|f^b\|_{L^2(S_z)}^2.
\]
Thus, integrating over $M$ and applying the dominated
convergence theorem, we obtain
\[
f^b=\sum_{m=0}^{\infty}P_m f^b
\quad\text{in }L^2(S),\qquad
\|f^b\|_{L^2(S)}^2
=\sum_{m=0}^{\infty}\|P_m f^b\|_{L^2(S)}^2.
\]
This establishes the desired orthogonal decomposition.
\end{proof}

Recall that $A^2(M,L^m)$ consists of holomorphic sections of $L^m$ which are square integrable with respect to $h^m$ and $dV_g$, as in \eqref{Eq Bergman inner product}. The following lemma identifies this Bergman space with $H_m^2(S)$.

\begin{lemma}\label{Lemma Hardy Bergman isomorphism}
For every integer $m\geq0$ and every $u\in H_m^2(S)$, let $f\in H^2(D)$ be the unique function with radial boundary value $u$. Then $f$ has the form $f(z,\xi)=b(z)\xi^m$ in the above trivialization over $U$, where $b$ is holomorphic on $U$. Moreover, the map
\[
I_m(u)=b(z)e_L^{m}\quad\text{over }U
\]
is globally well-defined and gives a unitary isomorphism $I_m:H_m^2(S)\longrightarrow A^2(M,L^m)$, where $e_L$ is the dual frame of $e_{L^*}$.
\end{lemma}

\begin{proof}
 First since $r_\theta(rv)=r\,r_\theta(v)$, for each fixed $\theta$ and almost every $v\in S$, we have
\[
(r_\theta^*f)^b(v)=\lim_{r\to1^-}f(r_\theta(rv))=\lim_{r\to1^-}f(r\,r_\theta(v))=f^b(r_\theta(v))=(r_\theta^*u)(v).
\]
Thus $r_\theta^*f$ and $e^{im\theta}f$ have the same boundary value $r_\theta^*u=e^{im\theta}u$. Uniqueness implies $r_\theta^*f=e^{im\theta}f$. Comparing coefficients in the local expansion of $f$ in \eqref{Eq Hardy fiber Taylor}, we obtain $f(z,\xi)=b(z)\xi^m$ for some holomorphic
function $b$ on $U$.

We next verify that $I_m$ is well-defined. Let $(\widetilde z,\widetilde\xi)$ be another trivialization, with $\widetilde z=F(z)$ and $\widetilde e_{L^*}=\eta(z)e_{L^*}$ on the overlap, where $\eta$ is holomorphic and nonvanishing. Then $\widetilde\xi=\eta(z)^{-1}\xi$ and $\widetilde e_L=\eta(z)^{-1}e_L$. Writing $f=\widetilde b(\widetilde z)\widetilde\xi^m=b(z)\xi^m$, we obtain $\widetilde b(F(z))\eta(z)^{-m}=b(z)$. Consequently, $\widetilde b(\widetilde z)\widetilde e_L^{m} =b(z)e_L^{m}$. Thus $I_m(u)$ is a globally defined holomorphic section of $L^m$.

 On $S \cap \pi^{-1}(U)$, we have $\xi=h(z)^{1/2}e^{i\theta}$ and hence $|u(z,\theta)|^2=|b(z)|^2h(z)^m=|I_m(u)(z)|_{h^m}^2$. We first integrate this equation in the angular variable $\theta$ and then over $M$, using a partition of unity subordinate to trivializing charts. This yields
\begin{equation}\label{Eq Hardy Bergman norm}
\|u\|_{L^2(S)}^2
=\int_M |I_m(u)(z)|_{h^m}^2dV_g(z)
=\|I_m(u)\|_{L^2(M,L^m)}^2.
\end{equation}
In particular, $I_m(u)\in A^2(M,L^m)$, and $I_m$ is a linear isometry.

It remains to prove the surjectivity of $I_m$. Given $s\in A^2(M,L^m)$, define
$\widehat{s}(v)=v^{\otimes m}(s(\pi(v)))$
for $v\in L^*$. When $m=0$, this means $\widehat s=s\circ\pi$. Locally, if $s=b(z)e_L^{m}$, then $\widehat s(z,\xi)=b(z)\xi^m$. Thus $\widehat s$ is holomorphic and
\[
\int_S|\widehat s(rv)|^2d\sigma(v)=r^{2m}\|s\|_{L^2(M,L^m)}^2.
\]
It follows that $\widehat{s}\in H^2(D)$.
Since $\widehat{s}(rv)=r^m\widehat{s}(v)$,
its radial boundary value is
$(\widehat{s})^b=\widehat{s}|_S\in H_m^2(S)$.
Thus $I_m((\widehat{s})^b)=s$, so $I_m$ is surjective. This completes the proof. 
\end{proof}

\begin{proof}[Proof of Theorem~\ref{Thm Szego Bergman relation}]
For each $m\geq0$, let $\{s_{m,\alpha}\}_\alpha$ be an orthonormal basis of $A^2(M,L^m)$. The index set may be finite or countably infinite, and is empty when this space is trivial. By Proposition~\ref{Prop Hardy Fourier decomposition} and Lemma~\ref{Lemma Hardy Bergman isomorphism}, the functions $\{\widehat s_{m,\alpha}|_S\}_{m,\alpha}$ form an orthonormal basis of $H^2(S)$. The orthonormal-basis formula in Proposition~\ref{Prop Szego basis expansion} implies
\begin{equation}\label{Eq Szego lifted basis}
\Pi(v,w)=\sum_{m=0}^{\infty}\sum_\alpha
\widehat s_{m,\alpha}(v)\overline{\widehat s_{m,\alpha}(w)},\qquad v,w\in D.
\end{equation}
In a local trivialization $(z, \xi)$ over $U\subset M$,  write $s_{m,\alpha}=b_{m,\alpha}(z)e_L^{m}$.
Then
\[
\mathcal B_m(z)=\sum_\alpha|b_{m,\alpha}(z)|^2h(z)^m,
\qquad
\widehat s_{m,\alpha}(z,\xi)=b_{m,\alpha}(z)\xi^m.
\]
As $1-\rho(v)=|\xi|^2h(z)^{-1}$, we obtain
\[
\sum_\alpha|\widehat s_{m,\alpha}(v)|^2
=\mathcal B_m(z)|\xi|^{2m}h(z)^{-m}
=\mathcal B_m(\pi(v))(1-\rho(v))^m.
\]
Setting $w=v$ in \eqref{Eq Szego lifted basis}, we obtain \eqref{Eq Szego Bergman series}.
Fix any compact subset $Z\Subset D$. The partial sums
\[
\sum_{m=0}^N \mathcal B_m(\pi(v))(1-\rho(v))^m
\]
are continuous and increase pointwise to the continuous function $\Pi(v,v)$ on $Z$. By Dini's theorem, the infinite series in \eqref{Eq Szego Bergman series} converges uniformly on $Z$.
\end{proof}

The same argument identifies the kernel of each Fourier component on the boundary. By Lemma~\ref{Lemma Hardy Bergman isomorphism}, the functions $\{\widehat{s}_{m,\alpha}|_S\}_\alpha$ form an orthonormal basis of $H_m^2(S)$. Let $\overline{\Pi}_m$ denote the orthogonal projection from $L^2(S)$ onto $H_m^2(S)$. Its kernel is
\[
\overline{\Pi}_m(x,y)
=\sum_\alpha
\widehat{s}_{m,\alpha}(x)
\overline{\widehat{s}_{m,\alpha}(y)},
\qquad x,y\in S.
\]
For fixed $m$, the above series converges locally
uniformly, together with all derivatives, and hence
defines a smooth kernel on $S\times S$.
On the diagonal, we have
\begin{equation}\label{Eq boundary Szego Bergman}
	\overline{\Pi}_m(x,x)=\mathcal B_m(\pi(x)),
	\qquad x\in S.
\end{equation}

\subsection{The asymptotic expansion of the Szeg\H{o} kernel}\label{Subsec Szego expansion}

In this subsection, we prove Theorem 1, following a similar argument as in the proof of \cite[Theorem 1.9]{EXX-BLF}. This approach is inspired by \cite{Ca99,EnZh}. By Theorem~\ref{Thm Szego Bergman relation},
\begin{equation}\label{Eq coefficient starting series}
\Pi(v,v)=\mathcal B_0(\pi(v))+\sum_{k=1}^{\infty}\mathcal B_k(\pi(v))(1-\rho(v))^k.
\end{equation}
All functions on $M$ appearing below are identified with their pullbacks to $D\cup S$ under $\pi$. Recall the following definition from \cite[Definition 3.1]{EXX-BLF}.
\begin{defn}\label{Def OM p}
	Let $F$ be a function on $\mathbb Z_+\times M$,
	smooth in the second variable. For a nonnegative integer $p$, we say $F\in\mathcal O_M(p)$ if, for every integer $l\geq0$ and every compact subset $X\Subset M$,
	there is a constant $C>0$ such that
	\[
	\|F(k,\cdot)\|_{C^l(X)}\leq Ck^{-p}, \quad \forall k\geq 1.
	\]
\end{defn}

We prove the following analogue of \cite[Proposition 3.2]{EXX-BLF}.
\begin{lemma}\label{Lemma Szego summed remainder}
If $F\in O_M(p)$ for an integer $p\geq2$, then
\[
\sum_{k=1}^{\infty}(1-\rho)^kF(k,\cdot)\in C^{p-2}(D\cup S).
\]
\end{lemma}

\begin{proof}
	Choose a local coordinate chart $(U,z)$ of $M$ and a holomorphic frame $e_{L^*}$ of $L^*$ over $U$, as in Section~\ref{Subsec Szego Bergman relation}. Let $\alpha=(\alpha_1,\alpha_2,\alpha_3,\alpha_4) \in\mathbb Z_{\geq 0}^n\times\mathbb Z_{\geq 0}^n \times\mathbb Z_{\geq 0}\times\mathbb Z_{\geq 0}$ be a multi-index with $|\alpha|\leq p-2$. Define the differential operator
	\[
	L_\alpha=
	\frac{\partial^{|\alpha|}}
	{\partial z^{\alpha_1}
		\partial\overline z^{\alpha_2}
		\partial\xi^{\alpha_3}
		\partial\overline\xi^{\alpha_4}}.
	\]
	Fix an open subset $V\Subset U$.
	Since $0\leq 1-\rho\leq 1$ on $D\cup S$ and $F\in\mathcal O_M(p)$, there exists some constant $C>0$, independent of $k$, such that
	\[
	\left|L_\alpha\bigl((1-\rho)^kF(k,\cdot)\bigr)\right|
	\leq Ck^{|\alpha|-p}\leq Ck^{-2}
	\]
	on $(D\cup S)\cap\pi^{-1}(V)$ for all $k\geq1$. Consequently, the series $\sum_{k=1}^\infty L_\alpha\bigl((1-\rho)^kF(k,\cdot)\bigr)$ converges absolutely and uniformly on	$(D\cup S)\cap\pi^{-1}(V)$
	for every $|\alpha|\leq p-2$. It follows that $\sum_{k=1}^\infty(1-\rho)^kF(k,\cdot)$ is of class $C^{p-2}$ on $D\cup S$. Since $U$ and $V$ are arbitrary, the proof is complete.	
\end{proof}

We next record some elementary identities for later use. For each nonnegative integer $p$, let $\lambda_{pi}$, $0\leq i\leq p$, be the real numbers defined by the polynomial identity in $x$:
\begin{equation}\label{Eq lambda definition}
x^p=\sum_{i=0}^p\lambda_{pi}\binom{x+i}{i}.
\end{equation}
Here $\binom{x+i}{i}=(x+1)\cdots(x+i)/i!$, with the convention that $\binom{x}{0}=1$. Since these polynomials form a basis, the $\lambda_{pi}$ are uniquely determined. Comparing leading coefficients in \eqref{Eq lambda definition}, we obtain $\lambda_{pp}=p!$. In addition, we set $\lambda_{pi}=0$ when $i>p$.

Similarly, for integers $0\leq i\leq j$, let $\tau_{ij}$ be the real numbers defined by
\begin{equation}\label{Eq tau definition}
-\frac{(\log(1-t))^i}{i!}=\sum_{j=i}^{\infty}\tau_{ij}t^j,\qquad |t|<1.
\end{equation}
In particular, $\tau_{00}=-1$, $\tau_{0j}=0$ for $j>0$, and
\begin{equation}\label{Eq tau diagonal}
\tau_{ii}=\frac{(-1)^{i+1}}{i!}.
\end{equation}

\begin{lemma}\label{Lemma Szego elementary sums}
For $0<t<1$ and every nonnegative integer $p$, we have
\begin{equation}\label{Eq polynomial sum}
\sum_{k=1}^{\infty}k^p(1-t)^k
=\sum_{i=0}^p\lambda_{pi}t^{-i-1}-\delta_{p0},
\end{equation}
where $\delta_{p0}$ is the Kronecker delta, i.e., it equals $1$ if $p=0$ and $0$ otherwise. For each nonnegative integer $i$, there exists a real analytic function $H_i$ near $t=0$ such that
\begin{equation}\label{Eq inverse power sum}
\sum_{k=1}^{\infty}\frac{(1-t)^k}{k^{i+1}}
=H_i(t)-\frac{(\log(1-t))^i}{i!}\log t.
\end{equation}
Consequently, for every integer $N\geq i$,
\begin{equation}\label{Eq inverse power finite reduction}
\sum_{k=1}^{\infty}\frac{(1-t)^k}{k^{i+1}}
-\left(\sum_{j=i}^N\tau_{ij}t^j\right)\log t
\end{equation}
extends as a $C^N$ function to $t=0$.
\end{lemma}

\begin{proof}

By taking the Taylor expansion of $(1-r)^{-i-1}$
at $r=0$, we have
\[
\sum_{k=0}^{\infty}\binom{k+i}{i}r^k
=(1-r)^{-i-1},\qquad |r|<1.
\]
If we substitute $t=1-r$, then
\[
\sum_{k=1}^{\infty}\binom{k+i}{i}(1-t)^k
=t^{-i-1}-1,\qquad 0<t<1.
\]
Combining this with \eqref{Eq lambda definition},
we obtain
\[
	\sum_{k=1}^{\infty}k^p(1-t)^k
	=\sum_{k=1}^{\infty}
	\sum_{i=0}^p\lambda_{pi}
	\binom{k+i}{i}(1-t)^k
	=\sum_{i=0}^p\lambda_{pi}
	(t^{-i-1}-1).
\]
Setting $x=0$ in \eqref{Eq lambda definition}
gives $\sum_{i=0}^p\lambda_{pi}=\delta_{p0}$.
This proves \eqref{Eq polynomial sum}.

To prove \eqref{Eq inverse power sum}, for $0<t<1$, denote
\begin{equation}\label{Eq Fi Ti definition}
	F_i(t)=\sum_{k=1}^{\infty}\frac{(1-t)^k}{k^{i+1}},
	\qquad
	T_i(t)=-\frac{(\log(1-t))^i}{i!}.
\end{equation}
We set $H_i(t)=F_i(t)-T_i(t)\log t$ for $0<t<1$, and show by induction on $i$ that $H_i$ extends as a real analytic function near $t=0$.

When $i=0$, the identities in \eqref{Eq Fi Ti definition} become $F_0(t)=-\log t$ and $T_0(t)=-1$,
thus $H_0=F_0-T_0\log t=0$. For $i\geq1$, we differentiate
\eqref{Eq Fi Ti definition} and obtain
\[
	F_i'(t)
	=-\sum_{k=1}^{\infty}
	\frac{(1-t)^{k-1}}{k^i}
	=-\frac{F_{i-1}(t)}{1-t},\qquad
	T_i'(t)
	=\frac{(\log(1-t))^{i-1}}{(i-1)!(1-t)}
	=-\frac{T_{i-1}(t)}{1-t}.
\]
Therefore,
\begin{equation}\label{Eq Hi recursive}
	\frac{d}{dt} H_i(t)=\frac{d}{dt}\big(F_i(t)-T_i(t)\log t\big) =-\frac{H_{i-1}(t)}{1-t}
	-\frac{T_i(t)}{t}.
\end{equation}
Assume $H_{i-1}$ is real analytic near $t=0$. Since $T_i(0)=0$ for $i\geq1$, the right-hand side of \eqref{Eq Hi recursive} is real analytic near $t=0$. Hence, integration shows that $H_i$ extends as a real analytic function near $t=0$. By induction, this holds for every nonnegative integer $i$.

Finally, we prove \eqref{Eq inverse power finite reduction}.
By \eqref{Eq tau definition}, we have
\[
T_i(t)=\sum_{j=i}^N\tau_{ij}t^j
+t^{N+1}R_{i,N}(t),
\]
where $R_{i,N}$ is real analytic near $t=0$.
Combining this with \eqref{Eq inverse power sum},
we obtain
\[
F_i(t)-\left(\sum_{j=i}^N\tau_{ij}t^j\right)
\log t
=H_i(t)+t^{N+1}R_{i,N}(t)\log t.
\]
Since $H_i$ and $R_{i,N}$ are real analytic near $t=0$ and $t^{N+1}\log t$ extends as a $C^N$ function to $t=0$, the proof is complete.
\end{proof}

\begin{proof}[Proof of Theorem \ref{Thm main expansion}]
\textit{Step 1. Expansions with remainders of finite regularity.} 
By \eqref{Eq Szego Bergman series}, we have
\begin{equation}\label{Eq Szego Bergman relation}
	\Pi(v,v)=\mathcal B_0(\pi(v))
	+\sum_{k=1}^{\infty}
	\mathcal B_k(\pi(v))(1-\rho(v))^k.
\end{equation}
Since $\mathcal B_0$ is smooth, to study the
singularity of $\Pi(v,v)$ along $S$, we only need to analyze the series over $k\geq1$ using the Tian--Yau--Zelditch--Catlin expansion \eqref{Eq general TYZ expansion}.

In view of Proposition \ref{Prop TYZ expansion}, for any integer $N\geq0$, we write
\begin{equation}\label{Eq truncated Bergman expansion}
	\mathcal B_k=\frac1{\pi^n}
	\sum_{\ell=0}^{n+N+1}a_\ell k^{n-\ell}
	+R_k^{(N)}, \qquad k\geq1.
\end{equation}

Substituting \eqref{Eq truncated Bergman expansion}
into \eqref{Eq Szego Bergman relation} and splitting
the finite sum at $\ell=n$, we obtain
\begin{equation}\label{Eq Szego finite sum splitting}
	\begin{aligned}
		\Pi(v,v)
		&=\mathcal B_0
		+\frac1{\pi^n}\sum_{k=1}^{\infty}(1-\rho)^k
		\sum_{\ell=0}^n a_\ell k^{n-\ell}
		+\frac1{\pi^n}\sum_{k=1}^{\infty}(1-\rho)^k
		\sum_{i=0}^N\frac{a_{n+1+i}}{k^{i+1}}
		+\sum_{k=1}^{\infty}(1-\rho)^kR_k^{(N)}\\
		&=:\mathcal B_0+(\mathrm I)+(\mathrm{II})
		+R^{(N)}.
	\end{aligned}
\end{equation}
With the notation in Definition~\ref{Def OM p}, the remainder estimates \eqref{Eq TYZ remainder estimates} give $R_k^{(N)}\in\mathcal O_M(N+2)$. By Lemma~\ref{Lemma Szego summed remainder}, we then have
\[
 R^{(N)}
=\sum_{k=1}^{\infty}(1-\rho)^kR_k^{(N)}
\in C^N(D\cup S).
\]

We first consider $(\mathrm I)$. By
\eqref{Eq polynomial sum}, we have
\[
\begin{aligned}
	(\mathrm I)
	=\frac1{\pi^n}\sum_{\ell=0}^n a_\ell
	\sum_{k=1}^{\infty}k^{n-\ell}(1-\rho)^k
	&=\frac1{\pi^n}\sum_{\ell=0}^n a_\ell
	\left(\sum_{i=0}^{n-\ell}
	\lambda_{n-\ell,i}\rho^{-i-1}
	-\delta_{n-\ell,0}\right)\\
	&=\rho^{-n-1}
	\left(\frac1{\pi^n}\sum_{\ell=0}^n a_\ell
	\sum_{i=0}^{n-\ell}
	\lambda_{n-\ell,i}\rho^{n-i}\right)
	-\frac{a_n}{\pi^n}.
\end{aligned}
\]
Define the polynomial in $t$ with smooth coefficients
in $z$ by
\begin{equation}\label{Eq P polynomial}
	P(z,t)=\frac1{\pi^n}\sum_{\ell=0}^n a_\ell(z)
	\sum_{i=0}^{n-\ell}\lambda_{n-\ell,i}t^{n-i}.
\end{equation}
Then
\[
(\mathrm I)
=\rho(v)^{-n-1}P(\pi(v),\rho(v))
-\frac{a_n(\pi(v))}{\pi^n}.
\]

We next consider $(\mathrm{II})$.
By \eqref{Eq inverse power finite reduction},
for every $0\leq i\leq N$, we have
\[
\sum_{k=1}^{\infty}
\frac{(1-\rho)^k}{k^{i+1}}
-\left(\sum_{j=i}^N\tau_{ij}\rho^j\right)
\log\rho\in C^N(D\cup S)
\]
Consequently,
\begin{multline*}
	(\mathrm{II})
	-\frac1{\pi^n}\sum_{i=0}^N a_{n+1+i}
	\left(\sum_{j=i}^N\tau_{ij}\rho^j\right)\log\rho\\
	=\frac1{\pi^n}\sum_{i=0}^N a_{n+1+i}
	\left(
	\sum_{k=1}^{\infty}\frac{(1-\rho)^k}{k^{i+1}}
	-\left(\sum_{j=i}^N\tau_{ij}\rho^j\right)
	\log\rho
	\right)\in C^N(D\cup S).
\end{multline*}
Interchanging the finite sums on the left-hand side, we obtain
\[
\begin{aligned}
	\frac1{\pi^n}\sum_{i=0}^N a_{n+1+i}
	\left(\sum_{j=i}^N\tau_{ij}\rho^j\right)
	=\sum_{j=0}^N
	\left(\frac1{\pi^n}\sum_{i=0}^j
	\tau_{ij}a_{n+1+i}\right)\rho^j.
\end{aligned}
\]
Define smooth functions on $M$ by
\begin{equation}\label{Eq beta coefficients}
	\beta_j(z)=\frac1{\pi^n}\sum_{i=0}^j
	\tau_{ij}a_{n+1+i}(z),\qquad j\geq0.
\end{equation}
Then
\[
(\mathrm{II})
-\left(\sum_{j=0}^N\beta_j\rho^j\right)
\log\rho\in C^N(D\cup S).
\]

Combining the above computations with \eqref{Eq Szego finite sum splitting}, we obtain
\begin{equation}\label{Eq theorem one finite expansion}
		\Pi(v,v)
		-\frac{P(\pi(v),\rho(v))}{\rho(v)^{n+1}}
		-\left(\sum_{j=0}^N
		\beta_j(\pi(v))\rho(v)^j\right)\log\rho(v)
		\in C^N(D\cup S).
\end{equation}

\textit{Step 2. Construction of smooth coefficients.}
Since $\{\beta_j\}_{j\geq0}$ is a sequence of
smooth functions on $M$, by Borel's theorem,
there exists a smooth function
$\Psi\in C^\infty(M\times\mathbb R)$ such that
\[
\partial_t^j\Psi(z,0)=j!\beta_j(z),
\qquad j\geq0.
\]
By the Taylor expansion of $\Psi(z,t)$ at $t=0$,
for every integer $N\geq0$, we have
\[
\Psi(z,t)=\sum_{j=0}^N\beta_j(z)t^j
+t^{N+1}E_N(z,t),
\]
where $E_N$ is smooth on $M\times\mathbb R$.

Define
\begin{equation}\label{Eq psi definition}
	\psi(v)=\Psi(\pi(v),\rho(v)),
	\qquad v\in D\cup S.
\end{equation}
Then $\psi\in C^\infty(D\cup S)$, and
\[
\left(\psi-\sum_{j=0}^N\beta_j\rho^j\right)
\log\rho
=\rho^{N+1}E_N(\pi(v),\rho(v))
\log\rho\in C^N(D\cup S).
\]
Combining this with
\eqref{Eq theorem one finite expansion},
we obtain that the function
\begin{equation}\label{Eq G definition}
	G(v):=\Pi(v,v)
	-\frac{P(\pi(v),\rho(v))}{\rho(v)^{n+1}}
	-\psi(v)\log\rho(v)
\end{equation}
belongs to $C^N(D\cup S)$ for every integer
$N\geq0$. Hence $G\in C^\infty(D\cup S)$.

Now define
\begin{equation}\label{Eq phi definition}
	\phi(v)=P(\pi(v),\rho(v))
	+\rho(v)^{n+1}G(v),\qquad v\in D\cup S.
\end{equation}
Since $P$ and $G$ are smooth on $D\cup S$, we have $\phi\in C^\infty(D\cup S)$. By \eqref{Eq G definition}, we obtain
\[
\Pi(v,v)=\frac{\phi(v)}{\rho(v)^{n+1}}
+\psi(v)\log\rho(v),\qquad v\in D.
\]
This proves the first part of Theorem~\ref{Thm main expansion}.

\textit{Step 3. Uniqueness of the coefficient functions.}
Suppose $\widetilde\phi,\widetilde\psi\in
C^\infty(D\cup S)$ also satisfy \eqref{Eq general Szego expansion}. Then
\[
(\phi-\widetilde\phi)
+\rho^{n+1}(\psi-\widetilde\psi)\log\rho=0
\qquad\text{on }D.
\]
Restricting this identity to smooth curves in $D$ intersecting $S$ transversely, as in \cite[Proof of Corollary 1.4]{EXX-BLF}, we conclude that both $\phi-\widetilde\phi$ and $\psi-\widetilde\psi$ vanish to infinite order along $S$. This proves the uniqueness part of Theorem~\ref{Thm main expansion}.
\end{proof}

\subsection{Explicit relations between the coefficients}\label{Subsec Szego coefficients}
In this subsection, we prove Theorem~\ref{Thm main coefficients} using the identities obtained in the construction of $\phi$ and $\psi$ in the proof of Theorem~\ref{Thm main expansion}. These identities give explicit relations between their Taylor coefficients and the Tian--Yau--Zelditch--Catlin coefficients.

\begin{proof}[Proof of Theorem~\ref{Thm main coefficients}]
	By \eqref{Eq phi definition} and \eqref{Eq psi definition} in the proof of Theorem~\ref{Thm main expansion}, we have
	\[
	\phi(v)=P(\pi(v),\rho(v))
	+\rho(v)^{n+1}G(v),
	\qquad
	\psi(v)=\Psi(\pi(v),\rho(v)),
	\]
	where $G\in C^\infty(D\cup S)$ and $\Psi\in C^\infty(M\times\mathbb R)$  satisfies $\partial_t^j\Psi(z,0)=j!\beta_j(z)$ for every $j\geq0$.
	
	Set
	\begin{equation}\label{Eq explicit A coefficients}
		A_j=\frac1{\pi^n}\sum_{\ell=0}^j
		\lambda_{n-\ell,n-j}a_\ell,
		\qquad 0\leq j\leq n.
	\end{equation}
	Then, by \eqref{Eq P polynomial}, we obtain
	\[
	P(z,t)=\frac1{\pi^n}\sum_{\ell=0}^n a_\ell(z)
	\sum_{i=0}^{n-\ell}\lambda_{n-\ell,i}t^{n-i}
	=\sum_{j=0}^n A_j(z)t^j.
	\]
	
	Moreover, by \eqref{Eq beta coefficients}, we have
	\begin{equation}\label{Eq explicit B coefficients}
		B_j:=\beta_j
		=\frac1{\pi^n}\sum_{i=0}^j
		\tau_{ij}a_{n+1+i},
		\qquad j\geq0.
	\end{equation}
	Consequently, for every integer $N\geq0$,
	\[
	\phi=\sum_{j=0}^n A_j\rho^j
	+O(\rho^{n+1}),
	\qquad
	\psi=\sum_{j=0}^N B_j\rho^j
	+O(\rho^{N+1}).
	\]
	The uniqueness of these coefficients follows
	from the uniqueness assertion in
	Theorem~\ref{Thm main expansion}.
	
	Finally, for every integer $m\geq0$,
	\eqref{Eq explicit B coefficients} expresses
	$(B_0,\ldots,B_m)^{\mathsf T}$ as a lower
	triangular linear transformation of
	$(a_{n+1},\ldots,a_{n+1+m})^{\mathsf T}$.
	By \eqref{Eq tau diagonal}, its diagonal
	entries are
	\[
	\frac{1}{\pi^n}\tau_{jj}=\frac{(-1)^{j+1}}{\pi^n j!},
	\qquad 0\leq j\leq m,
	\]
	which are all nonzero. Thus the transformation is invertible. Consequently, for any subset $U\subset M$, conditions (1) and (2) in Theorem~\ref{Thm main coefficients} are equivalent. This finishes the proof of Theorem~\ref{Thm main coefficients}.
\end{proof}

\begin{rmk}
The coefficients $A_j$ and $B_j$ can be computed
successively from \eqref{Eq explicit A coefficients}
and \eqref{Eq explicit B coefficients}.
The constants $\lambda_{pi}$ are determined
by the polynomial identity \eqref{Eq lambda definition},
while $\tau_{ij}$ are obtained from the Taylor
expansion \eqref{Eq tau definition}.
For example, since $\lambda_{n,n}=n!$ and $\tau_{00}=-1$, the leading coefficients are
\[
A_0=\frac{n!a_0}{\pi^n}=\frac{n!}{\pi^n},
\qquad
B_0=-\frac{a_{n+1}}{\pi^n}.
\]
\end{rmk}

\subsection{Proofs of the corollaries}\label{Sec Szego corollaries}

\begin{proof}[Proof of Corollary~\ref{Cor SLF criterion}]
By Theorem~\ref{Thm main coefficients}, $\psi$ vanishes to infinite order along $\pi^{-1}(U)\cap S$ if and only if $B_j|_U=0$ for every $j\geq0$, which is equivalent to $a_j|_U=0$ for every $j>n$. Taking $U=M$ proves the global assertion.
\end{proof}

We now prove Corollary~\ref{Cor constant coefficients}. Note that the metric in this corollary is not assumed to be globally polarized, so the Riemann--Roch argument in \cite[\S~4]{EXX-BLF} does not apply directly here. We instead use the decomposition of the Bergman coefficients obtained by Alexakis--Hirachi \cite{AH17}.

\begin{proof}[Proof of Corollary~\ref{Cor constant coefficients}]
By Definition~\ref{Def SLF metric}, (i) implies $a_j\equiv0$ for every $j>n$, so (i) implies (ii). Clearly, (ii) implies (iii). It remains to prove that (iii) implies (i). Suppose (iii) holds. By the decomposition of the Bergman coefficients obtained by Alexakis--Hirachi \cite[Main Theorem and \S~1.2]{AH17}, each $a_j$ differs from a normalized scalar contraction of the $(j,j)$-component of the Todd form of $T^{1,0}M$ by divergence terms. This decomposition is a universal local identity valid in every complex dimension. For $j>n$, the Todd component vanishes, so Stokes' theorem gives $\int_Ma_j\,dV_g=0$. Since each $a_j$ is continuous and has one sign on $M$, it must vanish identically. Definition~\ref{Def SLF metric} then implies (i).
\end{proof}

\begin{rmk}
One can also verify that the integrals of $\widetilde a_j$ depend only on the K\"ahler class, as required in the Main Theorem of Alexakis--Hirachi \cite{AH17}. Their universal decomposition and the same argument as above then remove the global polarization assumption in \cite[Corollary~1.13]{EXX-BLF}, with condition~{\rm(ii)} understood as Bergman logarithmic flatness of the metric $g$.
\end{rmk}

\begin{proof}[Proof of Corollary~\ref{Cor locally homogeneous}]
Recall that, by Proposition~\ref{Prop TYZ expansion}, the coefficients $a_j$ are polynomials in the curvature of $(M,g)$ and its covariant derivatives. Since $(M,g)$ is locally homogeneous, each $a_j$ is constant on $M$. Hence, by Corollary~\ref{Cor constant coefficients}, $g$ is Szeg\H{o} logarithmically flat. The sphericity assertion follows from the work of Webster \cite{Web78}, Kamishima \cite{Kam94}, Bryant \cite{Bry01} and Wang \cite[Theorem 12]{Wang}, as in \cite[Corollary 1.14]{EXX-BLF}. This finishes the proof of Corollary~\ref{Cor locally homogeneous}.
\end{proof}

\section{Bergman and Szeg\H{o} logarithmic flatness}\label{Sec BLF SLF}

In \S~\ref{Subsec KE comparison}, we prove Theorem \ref{Thm main comparison}. In \S~\ref{Subsec flatness examples}, we study Bryant's metrics to show that Bergman logarithmic flatness does not imply Szeg\H{o} logarithmic flatness in general and that the coefficient $\widetilde a_{n+1}$ (see its definition in \eqref{Eq two Bergman expansions}) need not vanish for a Bergman logarithmically flat metric.

\subsection{Comparison in the K\"ahler--Einstein case}\label{Subsec KE comparison}
Let $g$ be a K\"ahler metric on a complex manifold $M$ of dimension $n$, not necessarily globally polarized. On a sufficiently small relatively compact coordinate ball, choose a polarization $(L,h)$ of $g$. Write $\mathcal B_k$ and $\widetilde{\mathcal B}_k$ for the local Bergman kernel functions of $L^k$ and $L^k\otimes\mathcal{C}_M$ on the coordinate ball, respectively, where $\mathcal{C}_M$ denotes the canonical line bundle of $M$, equipped with the Hermitian metric induced by $g$. The $L^2$ inner products on sections of these bundles are defined using their induced Hermitian metrics and the volume form $dV_g$. Their expansions are denoted by
\begin{equation}\label{Eq two Bergman expansions}
\mathcal B_k\sim\frac1{\pi^n}\sum_{j=0}^\infty a_jk^{n-j},\qquad
\widetilde{\mathcal B}_k\sim\frac1{\pi^n}\sum_{j=0}^\infty\widetilde a_jk^{n-j},
\qquad a_0=\widetilde a_0=1.
\end{equation}
The coefficients $\widetilde a_j$ are those of the canonically twisted expansion in \cite[Equation (1.5)]{EXX-BLF}.  The coefficients $\widetilde a_j$ and $a_j$ are universal polynomials in the curvature of $g$ and its covariant derivatives.

Applying Theorem \ref{Thm main coefficients} and \cite[Theorem 1.12]{EXX-BLF} on relatively compact coordinate balls, we obtain, for every open set $U\subset M$,
\begin{equation}\label{Eq comparison criteria}
\begin{aligned}
g|_U\text{ is Szeg\H{o} logarithmically flat}&\quad\Longleftrightarrow\quad a_j|_U=0\quad(j\geq n+1),\\
g|_U\text{ is Bergman logarithmically flat}&\quad\Longleftrightarrow\quad \widetilde a_j|_U=0\quad(j\geq n+2).
\end{aligned}
\end{equation}
These criteria are purely local. When a global polarization satisfies the hypotheses of Proposition~\ref{Prop TYZ expansion}, the same coefficients occur in the corresponding global expansions.

In the following proposition, we use the generalized binomial coefficients
\[
\binom{r}{k}=\frac{r(r-1)\cdots(r-k+1)}{k!},\qquad r\in\mathbb Z,\quad k\geq1,
\]
with $\binom{r}{0}=1$. In particular, $\binom{r}{k}=0$ when $k>r\geq0$.

\begin{prop}\label{Prop KE shift}
Suppose that $g$ is K\"ahler--Einstein on $U$, with
\begin{equation}\label{Eq KE normalization}
-\partial\overline\partial\log\det g
=\lambda\,\partial\overline\partial\varphi,\qquad \lambda\in\mathbb R,
\end{equation}
where $h=e^{-\varphi}$. Then the coefficients $a_j$ and $\widetilde a_j$ in \eqref{Eq two Bergman expansions} satisfy
\begin{equation}\label{Eq KE shift}
a_j=\sum_{\ell=0}^j\binom{n-\ell}{j-\ell}\lambda^{j-\ell}\widetilde a_\ell,\qquad j\geq0.
\end{equation}
 In particular, $a_{n+1}=\widetilde a_{n+1}$ on $U$.
\end{prop}

\begin{proof}
Since the assertion is local, we may assume that $U$ is a relatively compact coordinate ball with coordinates $z=(z_1,\ldots,z_n)$ on $U$ and a holomorphic frame $e_L$ of $L$. Write $h(e_L,e_L)=e^{-\varphi}$.
By \eqref{Eq KE normalization}, $\log\det g+\lambda\varphi$ is pluriharmonic. Hence we can write
\begin{equation}\label{Eq KE determinant}
\det g=|F|^2e^{-\lambda\varphi},
\end{equation}
where $F$ is a nowhere-vanishing holomorphic function on $U$.
Write $dV_0$ for the Euclidean volume form and put $dz=dz_1\wedge\cdots\wedge dz_n$. In particular, $dV_g=\det g\,dV_0$ and $|dz|_{\mathcal C_M}^2=(\det g)^{-1}$.

For each real parameter $t>0$, consider the weighted Bergman spaces
\begin{equation}\label{Eq local weighted Bergman spaces}
\mathcal H_t=\left\{f\in\mathcal O(U):\int_U|f|^2e^{-t\varphi}\,dV_g<\infty\right\},\qquad
\widetilde{\mathcal H}_t=\left\{q\in\mathcal O(U):\int_U|q|^2e^{-t\varphi}\,dV_0<\infty\right\},
\end{equation}
with the corresponding weighted $L^2$ inner products. The corresponding Bergman kernel functions are given by
\[
\mathcal B_t(z)=\sum_\alpha|f_{t,\alpha}(z)|^2e^{-t\varphi(z)},\qquad
\widetilde{\mathcal B}_t(z)=\sum_\alpha|q_{t,\alpha}(z)|^2e^{-t\varphi(z)}(\det g(z))^{-1},
\]
where $\{f_{t,\alpha}\}_\alpha$ and $\{q_{t,\alpha}\}_\alpha$ are orthonormal bases of $\mathcal H_t$ and $\widetilde{\mathcal H}_t$, respectively. For positive integers $t=m$, the maps $f\mapsto f e_L^{m}$ and $q\mapsto q e_L^{m}\otimes dz$ identify these spaces unitarily with the Bergman spaces of $L^m$ and $L^m\otimes\mathcal C_M$ over $U$, respectively. Indeed, the factor $(\det g)^{-1}$ in the canonical metric cancels the factor $\det g$ in $dV_g$. The displayed kernel functions therefore agree with the corresponding bundle Bergman kernel functions. See \cite[Appendix]{EXX-BLF} for more details on these spaces and their kernels.

For every positive integer $m$ such that $m+\lambda>0$, multiplication by $F$ gives a unitary isomorphism $\mathcal H_m\to\widetilde{\mathcal H}_{m+\lambda}$, since
\[
\int_U|f|^2e^{-m\varphi}\,dV_g
=\int_U|Ff|^2e^{-(m+\lambda)\varphi}\,dV_0,
\]
and its inverse is multiplication by $F^{-1}$. Hence, if $\{f_\alpha\}_\alpha$ is an orthonormal basis of $\mathcal H_m$, then $\{Ff_\alpha\}_\alpha$ is an orthonormal basis of $\widetilde{\mathcal H}_{m+\lambda}$. Using this basis, we obtain
\[
\begin{aligned}
\widetilde{\mathcal B}_{m+\lambda}(z)
&=\sum_\alpha|F(z)f_\alpha(z)|^2e^{-(m+\lambda)\varphi(z)}(\det g(z))^{-1}\\
&=\sum_\alpha|f_\alpha(z)|^2e^{-m\varphi(z)}
=\mathcal B_m(z).
\end{aligned}
\]

The local construction and weighted $L^2$ estimates in \cite[Appendix]{EXX-BLF} apply to the weights $e^{-t\varphi}$ for all sufficiently large real $t$; they do not require $t$ to be an integer (see also \cite[Corollary 5.2]{HS22} which corresponds to the twisted case). 
Thus both expansions in \eqref{Eq two Bergman expansions} also hold with $k$ replaced by a real parameter $t$ tending to $+\infty$ (with the same universal coefficients). Taking $t=m+\lambda$ for the twisted kernel and $t=m$ for the untwisted kernel, we obtain for every $N\geq0$,
\[
\sum_{\ell=0}^N\widetilde a_\ell(m+\lambda)^{n-\ell}
=\sum_{j=0}^N a_jm^{n-j}
+O(m^{n-N-1}),
\]
locally uniformly on $U$. For $0\leq\ell\leq N$, Taylor expansion of $(1+x)^{n-\ell}$ at $x=0$, with $x=\lambda/m$, gives
\[
(m+\lambda)^{n-\ell}
=\sum_{j=\ell}^N\binom{n-\ell}{j-\ell}\lambda^{j-\ell}m^{n-j}
+O(m^{n-N-1})
\]
as $m\to+\infty$. Substituting and comparing coefficients yields
\[
a_j=\sum_{\ell=0}^j\binom{n-\ell}{j-\ell}\lambda^{j-\ell}\widetilde a_\ell,
\qquad 0\leq j\leq N.
\]
Since $N$ is arbitrary, this proves \eqref{Eq KE shift}. For $j=n+1$, all terms with $\ell\leq n$ vanish because $\binom{n-\ell}{n+1-\ell}=0$, leaving $a_{n+1}=\widetilde a_{n+1}$.
\end{proof}

\begin{thm}\label{Thm KE SLF BLF}
Let $g$ be a K\"ahler metric on $M$ that is K\"ahler--Einstein on an open set $U\subset M$. If $g|_U$ is Szeg\H{o} logarithmically flat, then $\widetilde a_j|_U=0$ for all $j\geq n+1$, and hence $g|_U$ is Bergman logarithmically flat. If $g|_U$ is Bergman logarithmically flat, then $g|_U$ is Szeg\H{o} logarithmically flat if and only if $\widetilde a_{n+1}|_U=0$.
\end{thm}

\begin{proof}
If $g|_U$ is Szeg\H{o} logarithmically flat, then \eqref{Eq comparison criteria} gives $a_j=0$ on $U$ for all $j\geq n+1$. By Proposition~\ref{Prop KE shift}, we have
\[
0=a_j=\sum_{\ell=n+1}^j\binom{n-\ell}{j-\ell}\lambda^{j-\ell}\widetilde a_\ell,
\qquad j\geq n+1,
\]
since $\binom{n-\ell}{j-\ell}=0$ for $\ell\leq n$ and $j\geq n+1$. The coefficient of $\widetilde a_j$ is $1$, so induction on $j$ gives $\widetilde a_j=0$ on $U$ for all $j\geq n+1$. Hence $g|_U$ is Bergman logarithmically flat by \eqref{Eq comparison criteria}.

Conversely, suppose that $g|_U$ is Bergman logarithmically flat. By \eqref{Eq comparison criteria}, $\widetilde a_j=0$ on $U$ for all $j\geq n+2$. Proposition~\ref{Prop KE shift} therefore gives
\[
a_j=\binom{-1}{j-n-1}\lambda^{j-n-1}\widetilde a_{n+1}
=(-\lambda)^{j-n-1}\widetilde a_{n+1},
\qquad j\geq n+1.
\]
Thus $a_j=0$ on $U$ for all $j\geq n+1$ if and only if $\widetilde a_{n+1}=0$ on $U$, proving the second assertion by \eqref{Eq comparison criteria}.
\end{proof}

\begin{proof}[Proof of Theorem~\ref{Thm main comparison}]
Apply Theorem~\ref{Thm KE SLF BLF} with $U=M$.
\end{proof}

\begin{rmk}\label{Rmk dimension one antecedents}
In complex dimension one, the comparison is related to the results of Hirachi \cite{Hi93a}, Hirachi--Komatsu--Nakazawa \cite{HKN93}, and Takeuchi \cite{Ta20}. In particular, \cite[Theorem 1.8]{Ta20} states that, for a relatively compact strictly pseudoconvex domain in a complex surface whose boundary admits a pseudo-Einstein contact form, the condition $\psi_\theta=O(\rho^3)$ is equivalent to the chosen contact form $\theta$ being pseudo-Einstein and the boundary being spherical. Here $\psi_\theta$ denotes the logarithmic coefficient in the diagonal expansion of the Szeg\H{o} kernel defined with respect to the volume form $\theta\wedge d\theta$, in analogy with \eqref{Eq general Szego expansion}. Hirachi \cite{Hi93a} had previously obtained this characterization for domains in $\mathbb C^2$ with transverse symmetry. For a compact polarized K\"ahler--Einstein Riemann surface, the connection form $\alpha$ introduced before \eqref{Eq circle measure} is a torsion-free contact form with constant Webster scalar curvature, hence is pseudo-Einstein. Moreover, $\alpha\wedge d\alpha$ is a constant multiple of $d\sigma$. Thus, taking $\theta=\alpha$, if $g$ is Szeg\H{o} logarithmically flat, then $\psi_\theta=O(\rho^3)$. By \cite[Theorem 1.8]{Ta20}, the circle bundle is spherical and hence Bergman logarithmically flat. The parameter-shift argument above applies to noncompact bases as well, in every dimension.
\end{rmk}

\subsection{Examples distinguishing the two notions of logarithmic flatness}\label{Subsec flatness examples}

The purpose of this subsection is twofold. First, we show that Bergman logarithmic flatness does not imply Szeg\H{o} logarithmic flatness. Second, we answer a question arising from our earlier work \cite{EXX-BLF}. Recall that in \cite{EXX-BLF}, we proved that every compact locally homogeneous polarized K\"ahler manifold is Bergman logarithmically flat. In fact, for these metrics we showed that $\widetilde a_j\equiv0$ for every $j\geq n+1$, whereas Bergman logarithmic flatness requires only that $\widetilde a_j\equiv0$ for $j\geq n+2$ (see the definition of $\widetilde a_j$ at the beginning of \S~\ref{Subsec KE comparison}). It is therefore natural to ask whether a Bergman logarithmically flat metric can have $\widetilde a_{n+1}\not\equiv0$. We answer this question affirmatively, even for metrics whose circle bundles are locally spherical. Both types of examples arise from Bryant's complete Bochner--K\"ahler metrics on $\mathbb C^n$ constructed in \cite{Bry01}.

We recall Bryant's family of complete Bochner--K\"ahler metrics on $\mathbb C^n$ from \cite[\S~4.4.5]{Bry01}.
Fix $n\geq2$ and $\mathbf p=(p_1,\ldots,p_n)\in[0,\infty)^n$. Following Bryant \cite[\S~4.4.5]{Bry01}, define $s=s_{\mathbf p}$ on $\mathbb C^n$ by
\begin{equation}\label{Eq Bryant potential}
s=\sum_{j=1}^n|z_j|^2e^{-p_js}.
\end{equation}
More precisely, the function $s=s_{\mathbf p}(z)$ is defined implicitly by \eqref{Eq Bryant potential} as follows: For each fixed $z$, the difference between the left-hand side and the right-hand side is strictly increasing in $s\geq0$, is nonpositive at $s=0$, and tends to infinity as $s\to\infty$. Hence there is a unique nonnegative solution, which depends real analytically on $z$ by the implicit function theorem. Moreover, by uniqueness, $s(z)$ depends only on $|z_1|,\ldots,|z_n|$.

Put
\begin{equation}\label{Eq Bryant notation}
u_j:=|z_j|^2e^{-p_js},\qquad
A:=1+\sum_jp_ju_j,\qquad
P:=\sum_jp_j,\qquad
Q:=\sum_jp_j^2u_j.
\end{equation}
In particular, $s=\sum_j u_j$ and $A\geq1$. 
Consider the metric $g=(g_{i\bar j}):=(\partial_i\partial_{\bar j}s)$. Differentiating \eqref{Eq Bryant potential} gives $s_{z_i}=e^{-p_is}\bar z_i/A$ and $s_{\bar z_j}=e^{-p_js}z_j/A$. Differentiating $s_{z_i}$ with respect to $\bar z_j$ then yields
\[
g_{i\bar j}=\frac{e^{-p_is}}A\delta_{ij}
-\frac{e^{-(p_i+p_j)s}\bar z_i z_j}{A^2}
\left(p_i+p_j-\frac QA\right).
\]
One can directly verify that the matrix inverse of $g$ is $G=(G_{ij})$, where
$$G_{ij}:=A\left(\delta_{ij}e^{p_is}+(p_i+p_j+p_ip_js)\bar z_i z_j\right).$$ This agrees with the inverse matrix in Bryant's construction in \cite[\S~4.4.5]{Bry01}. Thus $g$ coincides with Bryant's metric, which is complete and Bochner--K\"ahler on $\mathbb C^n$ by \cite[Theorem~10]{Bry01}.

We next compute its determinant. Let $D_0=\operatorname{diag}(e^{p_1s},\ldots,e^{p_ns})$, $V=(\bar z,\mathbf p\bar z)$, and $C=\left(\begin{smallmatrix}0&1\\1&s\end{smallmatrix}\right)$. Then $G=A(D_0+VCV^*)$ and
\[
V^*D_0^{-1}V=\begin{pmatrix}s&A-1\\A-1&Q\end{pmatrix},
\qquad \det(I_2+CV^*D_0^{-1}V)=A^2.
\]
The matrix determinant lemma gives $\det G=e^{Ps}A^{n+2}$, and therefore
\begin{equation}\label{Eq Bryant determinant}
\det g=e^{-Ps}A^{-(n+2)}.
\end{equation}

We will denote Bryant's metric on $\mathbb C^n$ interchangeably by $g$ and $g_{\mathbf p}$. Let $L=\mathbb C^n\times\mathbb C$ be the trivial line bundle over $\mathbb C^n$, equipped with the Hermitian metric $h=e^{-s}$. Then $g$ is the K\"ahler metric induced by $(L,h)$, with associated volume form $dV_g=\det g\, dV_0$. The dual disk and circle bundles of $(L,h)$ are
\begin{equation}\label{Eq Bryant bundles}
D_{\mathbf p}=\{(z,\xi):|\xi|^2e^{s(z)}<1\},\qquad
S_{\mathbf p}=\{(z,\xi):|\xi|^2e^{s(z)}=1\}.
\end{equation}
The defining function $\rho$ used above is given by $\rho(z,\xi)=1-|\xi|^2e^{s(z)}$.
The measure $d\sigma$ used above is given by $d\sigma=\frac{d\theta}{2\pi}dV_g$, where $\theta$ is the angular coordinate in the fiber.

\begin{prop}\label{Prop Bryant BLF}
For every $\mathbf p\in[0,\infty)^n$, $S_{\mathbf p}$ is  spherical. In particular, $g_{\mathbf p}$ is Bergman logarithmically flat, and $\widetilde a_j\equiv0$ for $j\geq n+2$.
\end{prop}

\begin{proof}
By \cite[Proposition~4]{Wang} (cf. \cite{Web78}), vanishing of the Bochner tensor of the base is equivalent to  sphericity of the polarized circle bundle when $n\geq2$. Since $g_{\mathbf p}$ is Bochner--K\"ahler by \cite[Theorem~10]{Bry01}, $S_{\mathbf p}$ is locally spherical and hence Bergman logarithmically flat. Definition~\ref{Def BLF metric} then implies that $g_{\mathbf p}$ is Bergman logarithmically flat. Since $g_{\mathbf p}$ is complete, \cite[Theorem~1.12]{EXX-BLF} gives $\widetilde a_j\equiv0$ for all $j\geq n+2$.
\end{proof}

In the following, $\mathcal B_m$ and $\widetilde{\mathcal B}_m$ always denote the Bergman kernel functions of $L^m$ and $L^m\otimes\mathcal C_{\mathbb C^n}$ over $\mathbb C^n$, respectively.
We first compute their values $\mathcal B_m(0)$ and $\widetilde{\mathcal B}_m(0)$ at the origin. When computing integrals over $\mathbb C^n$ below, we will often integrate first in the angular variables and then change from the squared radii $|z_j|^2$ to the variables $u_j$ defined in \eqref{Eq Bryant notation}. The inverse change is $|z_j|^2=u_j e^{p_j\sum_i u_i}$, with Jacobian
\[
\det\left(\frac{\partial |z_j|^2}{\partial u_k}\right)
=e^{Ps}\det\left(\delta_{jk}+p_ju_j\right)
=e^{Ps}A.
\]
Writing $r_j=|z_j|$, the Jacobian formula gives
\[
\prod_{j=1}^n d(r_j^2)=e^{Ps}A\,du,
\qquad du=du_1\cdots du_n.
\]
Let $F\geq0$ be measurable on $[0,\infty)^n$. Using polar coordinates, we first integrate in the angular variables and then make the change of variables above. Together with \eqref{Eq Bryant determinant}, this yields
\begin{equation}\label{Eq Bryant radial measures}
\begin{aligned}
\int_{\mathbb C^n}F(u(z))\,dV_0
&=\pi^n\int_{[0,\infty)^n}F(u)e^{Ps}A\,du,\\
\int_{\mathbb C^n}F(u(z))\,dV_g
&=\pi^n\int_{[0,\infty)^n}F(u)A^{-(n+1)}\,du.
\end{aligned}
\end{equation}

\begin{prop}\label{Prop Bryant exact kernels}
For every integer $m>P$,
\begin{equation}\label{Eq Bryant twisted exact}
\widetilde{\mathcal B}_m(0)=\frac{(m-P)^{n+1}}{\pi^n m}.
\end{equation}
For every integer $m\geq1$,
\begin{equation}\label{Eq Bryant untwisted exact}
\mathcal B_m(0)=\frac{m^n}{\pi^n F_{\mathbf p}(1/m)},\qquad~\text{where}~
F_{\mathbf p}(x)=\int_{[0,\infty)^n}e^{-\sum_jv_j}
\left(1+x\sum_jp_jv_j\right)^{-(n+1)}\,dv.
\end{equation}
\end{prop}

\begin{proof}
We first compute the integrals defining the squared norms of the constant sections. For $m>P$, applying \eqref{Eq Bryant radial measures} we have
\[
\begin{aligned}
\int_{\mathbb C^n}e^{-ms}\,dV_0
&=\pi^n\int_{[0,\infty)^n}e^{-(m-P)\sum_j u_j}
\left(1+\sum_jp_ju_j\right)\,du\\
&=\pi^n\left(\frac1{(m-P)^n}+\frac{P}{(m-P)^{n+1}}\right)
=\frac{\pi^n m}{(m-P)^{n+1}}.
\end{aligned}
\]
For $m\geq1$, applying \eqref{Eq Bryant radial measures} and then substituting $v_j=mu_j$, we obtain
\[
\begin{aligned}
\int_{\mathbb C^n}e^{-ms}\,dV_g
&=\pi^n\int_{[0,\infty)^n}e^{-m\sum_j u_j}
\left(1+\sum_jp_ju_j\right)^{-(n+1)}\,du\\
&=\pi^n m^{-n}F_{\mathbf p}(1/m).
\end{aligned}
\]
These integrals are finite in the given ranges of $m$. For the second integral, we note that $0<F_{\mathbf p}(1/m)\leq1$.

Let $f$ be an entire function with $f(0)=1$. As the weights $e^{-ms}$ and $e^{-ms}\det g$ depend only on $|z_1|,\ldots,|z_n|$, the mean-value inequality yields
\[
\int_{\mathbb C^n}|f|^2e^{-ms}\,dV_0
\geq\int_{\mathbb C^n}e^{-ms}\,dV_0,
\qquad
\int_{\mathbb C^n}|f|^2e^{-ms}\,dV_g
\geq\int_{\mathbb C^n}e^{-ms}\,dV_g.
\] Equality holds for $f\equiv1$, so the constant function $1$ minimizes both integrals. Since $s(0)=0$ and $\det g(0)=1$, the extremal characterization of the Bergman kernel yields
\begin{equation}\label{Eq Bryant reciprocal kernels}
\widetilde{\mathcal B}_m(0)
=\left(\int_{\mathbb C^n}e^{-ms}\,dV_0\right)^{-1},
\qquad
\mathcal B_m(0)
=\left(\int_{\mathbb C^n}e^{-ms}\,dV_g\right)^{-1}.
\end{equation}
The preceding calculations prove \eqref{Eq Bryant twisted exact} and \eqref{Eq Bryant untwisted exact}.
\end{proof}

Consider the subfamily $p_1=\cdots=p_n=a>0$, and denote the corresponding metric and circle bundle by $g_a$ and $S_a$, respectively.

\begin{prop}\label{Prop Bryant equal parameters}
The Ricci eigenvalues relative to $g_a$ are bounded below by $-2a$. Consequently, $g_a$ satisfies the hypotheses of Proposition~\ref{Prop TYZ expansion}.
\end{prop}

\begin{proof}
Put $t=|z|^2$. Then $t=se^{as}$ and $A=1+as$ by \eqref{Eq Bryant notation}. For $z\ne0$, the tangential and radial eigenvalues of $(g_{i\bar j})$ in Euclidean coordinates are, respectively,
\[
s'(t)=\frac{e^{-as}}A,\qquad
s'(t)+ts''(t)=\frac{e^{-as}}{A^3}.
\]
Set $f(s)=\log\det g_a=-nas-(n+2)\log A$. Differentiating $f(s(t))$, we find that the eigenvalues of the Hermitian matrix $(-\partial_i\partial_{\bar j}\log\det g_a)$ relative to $(g_{i\bar j})$ are
\begin{equation}\label{Eq Bryant Ricci eigenvalues}
\begin{aligned}
\mu_{\mathrm{tan}}&=-f_s=na+\frac{(n+2)a}A,\\
\mu_{\mathrm{rad}}&=-f_s-sAf_{ss}=-2a+\frac{2(n+2)a}A.
\end{aligned}
\end{equation}
At the origin, all these eigenvalues equal $2(n+1)a$. Since $A\geq1$, they are bounded below by $-2a$ everywhere. With our normalization of the Ricci and K\"ahler forms, this gives
\[
\Ric(g_a)\geq-2a\,\omega_{g_a}.
\]
Since $g_a$ is complete, both hypotheses of Proposition~\ref{Prop TYZ expansion} hold.
\end{proof}

Recall that $a_j$ and $\widetilde a_j$ denote the local untwisted and canonically twisted Bergman coefficients of $g_{\mathbf p}$, respectively, as in \eqref{Eq two Bergman expansions}. They are universal polynomials in the curvature of $g_{\mathbf p}$ and its covariant derivatives. Their relation to the asymptotic expansions of the global kernel functions $\mathcal B_m$ and $\widetilde{\mathcal B}_m$ will be discussed below.

\begin{prop}\label{Prop Bryant sharp truncation}
For the Bryant metric $g_{\mathbf p}$  with $\mathbf p\ne0$, we have $\widetilde a_j\equiv0$ for $j\geq n+2$, but $\widetilde a_{n+1}(0)=(-1)^{n+1}P^{n+1}\ne0$. Thus the coefficient of $m^{-1}$ in the canonically twisted Bergman kernel expansion need not vanish for a Bergman logarithmically flat metric, even when its circle bundle is locally spherical.
\end{prop}

\begin{proof}
The vanishing of $\widetilde a_j$ for $j\geq n+2$ follows from Proposition~\ref{Prop Bryant BLF}. By \eqref{Eq Bryant twisted exact}, for all integers $m>P$,
\[
\pi^n\widetilde{\mathcal B}_m(0)
=\frac{(m-P)^{n+1}}m
=\sum_{j=0}^{n+1}\binom{n+1}{j}(-P)^jm^{n-j}.
\]
Since $g$ is complete, $\widetilde{\mathcal B}_m$ has an asymptotic expansion with coefficients $\widetilde a_j$ by \cite[Proposition~1.5 and Remark~1.7]{EXX-BLF}. Comparing coefficients in the preceding identity gives $\widetilde a_j(0)=\binom{n+1}{j}(-P)^j$ for $0\leq j\leq n+1$. Since $\mathbf p\ne0$ and all $p_j\geq0$, we have $P>0$ and hence
$\widetilde a_{n+1}(0)=(-1)^{n+1}P^{n+1}\ne0$.
\end{proof}

The coefficient $\widetilde a_{n+1}$ also determines whether products of Bergman logarithmically flat metrics remain Bergman logarithmically flat; see Theorem~\ref{Thm BLF product flatness}. We next show that Bryant's metrics with $\mathbf p\ne0$ are not Szeg\H{o} logarithmically flat when $n\geq2$. To this end, we first establish the following consequence of \eqref{Eq Bryant untwisted exact}.

\begin{lemma}\label{Lemma Bryant nontermination}
If $\mathbf p\ne0$, the asymptotic expansion of $\mathcal B_m(0)$ contains infinitely many nonzero terms with negative powers of $m$.
\end{lemma}

\begin{proof}
Set $Z(v)=\sum_jp_jv_j$. For $x\geq0$, the absolute value of the $k$th derivative with respect to $x$ of the integrand in \eqref{Eq Bryant untwisted exact} is bounded by $\frac{(n+k)!}{n!}e^{-\sum_jv_j}Z(v)^k$, which is integrable on $[0,\infty)^n$. We may therefore differentiate under the integral sign to any order. Applying Taylor's theorem to $F_{\mathbf p}$ at $x=0$, we obtain, for  $N\geq0$,
\[
F_{\mathbf p}(x)=\sum_{k=0}^N f_kx^k+O(x^{N+1}),
\qquad x\to0^+,
\]
where
\[
f_k=(-1)^k\binom{n+k}{k}
\int_{[0,\infty)^n}e^{-\sum_jv_j}Z(v)^k\,dv.
\]
Set $p_*:=\max_jp_j>0$ and
choose $j_0$ with $p_{j_0}=p_*$. Since $Z(v)\geq p_*v_{j_0}$, we have
\[
|f_k|\geq\binom{n+k}{k}(p_*)^k
\int_0^\infty e^{-t}t^k\,dt
=\frac{(n+k)!}{n!}(p_*)^k.
\]
Consequently, the Taylor series $\sum_{k\geq0}f_kx^k$ has radius of convergence zero.

Since $F_{\mathbf p}(0)=1$, the asymptotic series of $1/F_{\mathbf p}(x)$ is the formal inverse of this series. If it had only finitely many nonzero coefficients, it would be a polynomial $Q_0(x)$ with $Q_0(0)=1$. The series $\sum_{k\geq0}f_kx^k$ would then be the Taylor series of $1/Q_0(x)$, which is analytic near $0$. This contradicts its zero radius of convergence. Hence the expansion of $1/F_{\mathbf p}(x)$ has infinitely many nonzero coefficients. By \eqref{Eq Bryant untwisted exact}, the expansion of $\mathcal B_m(0)$ contains infinitely many nonzero terms with negative powers of $m$.
\end{proof}

\begin{prop}\label{Prop Bryant not SLF}
Let $n\geq2$ and $\mathbf p\ne0$. Then $g_{\mathbf p}$ is Bergman logarithmically flat but is not Szeg\H{o} logarithmically flat. 
\end{prop}

\begin{proof}
Bergman logarithmic flatness follows from Proposition~\ref{Prop Bryant BLF}. It remains to show that $g_{\mathbf p}$ is not Szeg\H{o} logarithmically flat. Lemma~\ref{Lemma Bryant nontermination} gives an asymptotic expansion
\[
\mathcal B_m(0)\sim\frac1{\pi^n}\sum_{j=0}^\infty c_jm^{n-j},
\]
where $c_j \in \mathbb C$ and $c_j\ne0$ for infinitely many $j>n$. We first show that these coefficients satisfy $c_j=a_j(0)$ for every $j\geq0$. Let $B$ be a ball centered at $0$, and let $\mathcal B_m^B$ denote the corresponding Bergman kernel function over $B$. The mean-value argument in the proof of Proposition~\ref{Prop Bryant exact kernels} gives
\begin{equation}\label{Eq Bryant ball kernel}
\mathcal B_m^B(0)=\left(\int_B e^{-ms}\,dV_g\right)^{-1}.
\end{equation}
By the defining equation for $s$ in \eqref{Eq Bryant potential}, $s(z)\to+\infty$ as $|z|\to\infty$, and $s$ vanishes only at $0$. Hence $s\geq\delta$ outside $B$ for some $\delta>0$. For $m\geq1$,
\begin{equation}\label{Eq Bryant tail estimate}
0\leq\int_{\mathbb C^n\setminus B}e^{-ms}\,dV_g
\leq e^{-(m-1)\delta}\int_{\mathbb C^n}e^{-s}\,dV_g,
\end{equation}
where the last integral is finite by the calculation in the proof of Proposition~\ref{Prop Bryant exact kernels}. That calculation also gives
\[
\int_{\mathbb C^n}e^{-ms}\,dV_g
=\pi^n m^{-n}F_{\mathbf p}(1/m)
\sim\pi^n m^{-n},
\]
since $F_{\mathbf p}(1/m)\to1$ as $m\to+\infty$. By \eqref{Eq Bryant tail estimate}, the same asymptotic holds for the integral over $B$. Thus, by \eqref{Eq Bryant ball kernel} and \eqref{Eq Bryant reciprocal kernels}, we have
\[
\begin{aligned}
\mathcal B_m^B(0)-\mathcal B_m(0)
&=\frac{\displaystyle\int_{\mathbb C^n\setminus B}e^{-ms}\,dV_g}
{\displaystyle\left(\int_B e^{-ms}\,dV_g\right)
\left(\int_{\mathbb C^n}e^{-ms}\,dV_g\right)}\\
&=O\!\left(m^{2n}e^{-(m-1)\delta}\right)
=O(m^{-\infty}).
\end{aligned}
\]
On $B$, the Ricci curvature of $g_{\mathbf p}$ is bounded below, and $B$ admits a complete K\"ahler metric. Thus Proposition~\ref{Prop TYZ expansion} applies, and uniqueness of the asymptotic expansion gives $c_j=a_j(0)$. Lemma~\ref{Lemma Bryant nontermination} and \eqref{Eq comparison criteria} now imply that $g_{\mathbf p}$ is not Szeg\H{o} logarithmically flat.
\end{proof}

We illustrate an explicit computation in the case $n=2$ and $\mathbf p=(a,b)$ with $a,b\geq0$. Using the formula for $f_k$ in the proof of Lemma~\ref{Lemma Bryant nontermination}, we obtain $$F_{(a,b)}(x)=1-3(a+b)x+12(a^2+ab+b^2)x^2-60(a^3+a^2b+ab^2+b^3)x^3+O(x^4)$$ as $x\to0^+.$ Inverting this series and using \eqref{Eq Bryant untwisted exact}, we get
\begin{equation}\label{Eq Bryant dimension two}
\mathcal B_m(0)=\frac1{\pi^2}\left[m^2+3(a+b)m-3(a-b)^2
+\frac{c(a,b)}m+O(m^{-2})\right],
\end{equation}
where
\[
c(a,b)=15(a^3+b^3)-3ab(a+b)
=12(a^3+b^3)+3(a+b)(a-b)^2.
\]
By the coefficient identification in the preceding proof, $a_3(0)=c(a,b)>0$ whenever $(a,b)\ne(0,0)$.
For $a=b>0$, Proposition~\ref{Prop Bryant equal parameters} ensures that Theorem~\ref{Thm main coefficients} applies, giving
\[
\psi(0,e^{i\theta})=-\frac{a_3(0)}{\pi^2}
=-\frac{24a^3}{\pi^2}\ne0.
\]
For this same example, \eqref{Eq Bryant twisted exact} gives
\[
\widetilde{\mathcal B}_m(0)
=\frac1{\pi^2}\left(m^2-6am+12a^2-\frac{8a^3}{m}\right),
\qquad m>2a.
\]
Hence $g_a$ is Bergman logarithmically flat but not Szeg\H{o} logarithmically flat, with $\widetilde a_3(0)=-8a^3\ne0$.

We finally prove Theorem~\ref{Thm Bryant_example}.

\begin{proof}[Proof of Theorem~\ref{Thm Bryant_example}]
Fix $a>0$ and take $g_a=g_{(a,\ldots,a)}$ on $\mathbb C^n, n \geq 2$. The metric is complete by \cite[Theorem~10]{Bry01}. By Proposition~\ref{Prop Bryant equal parameters}, it satisfies condition (b) in Proposition~\ref{Prop TYZ expansion}. The associated circle bundle $S_a$ is spherical by Proposition~\ref{Prop Bryant BLF}. Finally, by Proposition~\ref{Prop Bryant not SLF},
$g_a$ is Bergman logarithmically flat but not Szeg\H{o} logarithmically flat. 
\end{proof}

\begin{rmk}
A local CR equivalence between $S_{\mathbf p}$ and the sphere $\mathbb S^{2n+1}\subset\mathbb C^{n+1}$ need not take $d\sigma$ to the standard surface measure on the sphere. Therefore, sphericity alone does not imply the vanishing of the Szeg\H{o} logarithmic coefficient. The independence of Szeg\H{o} logarithmic flatness from the choice of local polarization follows from the curvature formulas for $a_j$, as explained in Remark~\ref{Rmk polarization independence}.
For $\mathbf p=0$, the metric $g_0$ is Euclidean. By the Bargmann--Fock formula, $\mathcal B_m=m^n/\pi^n$ for $m\geq1$, and $\mathcal B_0=0$. Hence $g_0$ is Szeg\H{o} logarithmically flat, and its diagonal Szeg\H{o} kernel is a rational function of $1-\rho$.
\end{rmk}

\subsection{Two conjectures}

Corollary \ref{Cor constant coefficients} motivates the following conjecture for noncompact homogeneous metrics.

\begin{conj}\label{Conj homogeneous SLF}
Let $(M,g;L,h)$ be a polarized manifold. If $(M,g)$ is homogeneous, then $g$ is Szeg\H{o} logarithmically flat.
\end{conj}

For a homogeneous polarized metric, the hypotheses of Proposition~\ref{Prop TYZ expansion} hold. Thus, by Corollary~\ref{Cor SLF criterion}, the conjecture says that the logarithmic term of the corresponding Szeg\H{o} kernel vanishes to infinite order along $S$.

Motivated by Theorem~\ref{Thm main comparison} and Takeuchi's theorem \cite[Theorem 1.8]{Ta20}, we propose the following conjecture for strictly pseudoconvex boundaries admitting pseudo-Einstein contact forms.

\begin{conj}\label{Conj pseudo Einstein comparison}
Let $\Omega$ be a relatively compact domain with smooth strictly pseudoconvex boundary in a complex manifold of dimension $n+1$, where $n\geq 2$. Suppose that $\partial\Omega$ admits a pseudo-Einstein contact form $\theta$. If the Szeg\H{o} kernel defined with respect to $\theta\wedge(d\theta)^n$ has no logarithmic singularity, then $\partial\Omega$ is Bergman logarithmically flat.
\end{conj}

When $n=1$, the corresponding statement follows from
\cite[Theorem 1.8]{Ta20}, as discussed in Remark~\ref{Rmk dimension one antecedents}.

\section{Product metrics and locally Hermitian symmetric metrics}\label{Sec products symmetric}
In this section, we prove Theorems~\ref{Thm SLF product flatness} and~\ref{Thm BLF product flatness}, as well as
Corollary~\ref{Cor Hermitian symmetric flatness}. Throughout this section, $a_j(g)$ and $\widetilde a_j(g)$ denote the coefficients in \eqref{Eq two Bergman expansions}. They are universal polynomials in the curvature of $g$ and its covariant derivatives. We use the
local characterizations in Definitions~\ref{Def SLF metric} and~\ref{Def BLF metric},
equivalently \eqref{Eq comparison criteria}.

\subsection{Rescaling and product formulas}\label{Subsec scaling products}
We prove two lemmas in this subsection.
\begin{lemma}\label{Lemma symmetric scaling}
	Let $g$ be a K\"ahler metric on a complex manifold $M$. For every constant $c>0$,
	\begin{equation}\label{Eq symmetric scaling}
		a_j(cg)=c^{-j}a_j(g),\qquad
		\widetilde a_j(cg)=c^{-j}\widetilde a_j(g),\qquad j\geq0.
	\end{equation}
\end{lemma}
\begin{proof}
	Write $n=\dim_{\mathbb C}M$. Since the coefficients are local invariants, it suffices to work
	on a sufficiently small coordinate ball $U\subset M$. Choose a K\"ahler potential $\varphi
	$ for $g$ on a neighborhood of $\overline U$, and use the trivial polarization with metric
	$h=e^{-\varphi}$. Let $dV_0$ denote Euclidean volume in these coordinates, and write $dz=dz
	_1\wedge\cdots\wedge dz_n$. Then the induced metric on the canonical bundle $\mathcal C_M$ by $g$ satisfies $|dz|_{\mathcal C_M
	}^2=(\det g)^{-1}$, and $dV_g=(\det g)dV_0$.
	For a real weight parameter $t>0$, let $\mathcal H_t^g$ and $\widetilde{\mathcal H}_t^g$ denote the weighted Bergman spaces defined in \eqref{Eq local weighted Bergman spaces}, where we include the superscript $g$ to indicate their dependence on the metric. Choose orthonormal bases $\{f_{t,\alpha}\}$ and $\{q_{t,\alpha}\}$ for $\mathcal H_t^g$ and $\widetilde{\mathcal H}_t^g$, respectively. Recall that for positive integers $t=m$, the maps $f\mapsto f e_L^m$ and $q\mapsto q e_L^m\otimes dz$ identify these spaces unitarily with the Bergman spaces of $L^m$ and $L^m\otimes\mathcal C_M$ over $U$, respectively. For the second identification, the factor $(\det g)^{-1}$ in the canonical metric cancels the factor $\det g$ in $dV_g$. We write the corresponding kernel functions as $\mathcal B_t^g$ and $\widetilde{\mathcal B}_t^g$ to indicate their dependence on the metric. Then they are given by
    \begin{equation}\label{Eq local kernel bases}
        \mathcal B_t^g=\sum_\alpha|f_{t,\alpha}|^2e^{-t\varphi},\qquad \widetilde{\mathcal B}_t^g=\sum_\alpha|q_{t,\alpha}|^2e^{-t\varphi}(\det g)^{-1}.
    \end{equation}
	For $cg$, the potential is $c\varphi$, the volume form is $c^n dV_g$, and the metric on the canonical bundle $\mathcal C_M$ is multiplied by $c^{-n}$. Hence $\{c^{-n/2}f_{cm,\alpha}\}$ is an orthonormal basis of $\mathcal H_m^{cg}$. Since $\widetilde{\mathcal H}_m^{cg}$ and $\widetilde{\mathcal H}_{cm}^g$ are both defined using the weight $e^{-cm\varphi}\,dV_0$, their inner products agree.  Hence $\{q_{cm,\alpha}\}$ is also an orthonormal basis of $\widetilde{\mathcal H}_m^{cg}$. Using these bases in \eqref{Eq local kernel bases}, together with $\det(cg)=c^n\det g$, we obtain
\begin{equation}\label{Eq symmetric exact scaling}
\mathcal B_m^{cg}=c^{-n}\mathcal B_{cm}^g,
\qquad
\widetilde{\mathcal B}_m^{cg}=c^{-n}\widetilde{\mathcal B}_{cm}^g.
\end{equation}
As explained in the proof of Proposition~\ref{Prop KE shift}, the local weighted expansions hold for real $t\to+\infty$. Taking $t=cm$, we obtain
	\[
	c^{-n}\mathcal B_{cm}^g
	\sim\frac1{\pi^n}\sum_{j=0}^\infty c^{-j}a_j(g)m^{n-j},
	\qquad
	c^{-n}\widetilde{\mathcal B}_{cm}^g
	\sim\frac1{\pi^n}\sum_{j=0}^\infty c^{-j}\widetilde a_j(g)m^{n-j}.
	\]
	Comparing these expansions with those for the left-hand sides of \eqref{Eq symmetric exact
		scaling} yields \eqref{Eq symmetric scaling}.
	\end{proof}

	\begin{lemma}\label{Lemma symmetric products}
	Let $(M_i,g_i)$ be K\"ahler manifolds, $i=1,2$, and let $g=g_1\oplus g_2$ on $M_1\times M_2$.
	Then, for every $j\geq0$,
	\begin{equation}\label{Eq symmetric products}
		a_j(g)(x_1,x_2)=\sum_{\ell=0}^j a_\ell(g_1)(x_1)a_{j-\ell}(g_2)(x_2),\qquad \widetilde a_j(g)
		(x_1,x_2)=\sum_{\ell=0}^j\widetilde a_\ell(g_1)(x_1)\widetilde a_{j-\ell}(g_2)(x_2).
	\end{equation}
\end{lemma}

\begin{proof}
Write $n_i=\dim_{\mathbb C}M_i$, $i=1,2$, and put $n=n_1+n_2=\dim_{\mathbb C}(M_1\times M_2)$. Choose sufficiently small coordinate balls $U_i\subset M_i$ and local polarizations $(L_i,h_i)$ defined on neighborhoods of $\overline{U_i}$. Their tensor product then polarizes $g=g_1\oplus g_2$ on $U_1\times U_2$.
The volume form is $dV_g=dV_{g_1}\,dV_{g_2}$. The canonical bundle of $U_1\times U_2$ is the tensor product of the canonical bundles of the two factors, and its induced metric is the product of their induced metrics.
Thus the untwisted and canonically twisted Bergman spaces on $U_1\times U_2$ are the Hilbert tensor products of the corresponding spaces on the factors. Products of orthonormal basis sections consequently form orthonormal bases on $U_1\times U_2$. This yields
	\begin{equation}\label{Eq symmetric exact products}
	 	\mathcal B_m^g(x_1,x_2)=\mathcal B_m^{g_1}(x_1)\mathcal B_m^{g_2}(x_2),\qquad \widetilde{\mathcal B}_m^g(x_1,x_2)=\widetilde{\mathcal B}_m^{g_1}(x_1)\widetilde{\mathcal B}_m^{g_2}(x_2).
	\end{equation}
In the first identity of \eqref{Eq symmetric exact products}, the two factors on the right have the local asymptotic expansions
\[
	\mathcal B_m^{g_1}(x_1)\sim\frac1{\pi^{n_1}}\sum_{\ell=0}^{\infty}a_\ell(g_1)(x_1)m^{n_1-\ell},\qquad \mathcal B_m^{g_2}(x_2)\sim\frac1{\pi^{n_2}}\sum_{r=0}^{\infty}a_r(g_2)(x_2)m^{n_2-r}.
	\]
Multiplying these expansions and comparing the coefficients of $m^{n-j}$ in the first identity of \eqref{Eq symmetric exact products} gives the first formula in \eqref{Eq symmetric products}. Applying the same argument to the twisted expansions gives the second formula.
\end{proof}

\subsection{Products of logarithmically flat metrics}\label{Subsec flatness products}
	Before proving Theorems~\ref{Thm SLF product flatness} and~\ref{Thm BLF product flatness}, we give the following remark about Bergman logarithmic flatness.

\begin{rmk}\label{Rmk BLF product counterexample}
Bergman logarithmic flatness is not preserved under products. Consider again the metrics $g_a$ and circle bundles $S_a$ introduced in \S~\ref{Subsec flatness examples}, corresponding to $p_1=\cdots=p_n=a>0$, with $n\geq2$. By Proposition~\ref{Prop Bryant BLF}, $S_a$ is locally spherical and $g_a$ is Bergman logarithmically flat.

Take $n=2$. By \eqref{Eq Bryant twisted exact}, the twisted kernel, evaluated at the origin, is
	\[
	\widetilde{\mathcal B}_m^{g_a}(0)=\frac{(m-2a)^3}{\pi^2 m}
	=\frac1{\pi^2}\left(m^2-6am+12a^2-\frac{8a^3}{m}\right),\qquad m>2a.
	\]
In particular,  $\widetilde a_3(g_a)(0)=-8a^3$. The identity \eqref{Eq product first obstruction}, proved below, gives $\widetilde a_6(g_a\oplus g_a)(0,0)=64a^6\neq0$. Since the product has complex dimension $4$, the metric $g_a\oplus g_a$ is not Bergman logarithmically flat.
\end{rmk}

\begin{proof}[Proof of Theorem~\ref{Thm SLF product flatness}]
Write $n_i=\dim_{\mathbb C}M_i$, $i=1,2$, and put $n=n_1+n_2$. Since $g_i$ is Szeg\H{o} logarithmically flat,  by \eqref{Eq comparison criteria}, $a_j(g_i)\equiv0$ for $j>n_i$, $i=1,2$. 
Therefore, a summand in the first formula of \eqref{Eq symmetric products} can be nonzero only if $\ell\leq n_1$ and $j-\ell\leq n_2$. It follows that $a_j(g_1\oplus g_2)\equiv0$ for $j>n$.
The criterion for Szeg\H{o} logarithmic flatness in \eqref{Eq comparison criteria} completes the proof.
\end{proof}

\begin{proof}[Proof of Theorem~\ref{Thm BLF product flatness}]
Put $n=n_1+n_2$. Since each $g_i$ is Bergman logarithmically flat, by \eqref{Eq comparison criteria}, $\widetilde a_j(g_i)\equiv0$ for $j\geq n_i+2$. Therefore, a summand in the second formula of \eqref{Eq symmetric products} can be nonzero only if $\ell\leq n_1+1$ and $j-\ell\leq n_2+1$. It follows that $\widetilde a_j(g_1\oplus g_2)\equiv0$ for $j>n+2$. At $j=n+2$, the only possibly nonzero summand gives
\begin{equation}\label{Eq product first obstruction}
\widetilde a_{n+2}(g_1\oplus g_2)(x_1,x_2)=\widetilde a_{n_1+1}(g_1)(x_1)\widetilde a_{n_2+1}(g_2)(x_2).
\end{equation}
By \eqref{Eq comparison criteria}, $g_1\oplus g_2$ is therefore Bergman logarithmically flat if and only if
\begin{equation}\label{Eq BLF product condition}
\widetilde a_{n_1+1}(g_1)(x_1)\widetilde a_{n_2+1}(g_2)(x_2)=0 \quad\text{for every }(x_1,x_2)\in M_1\times M_2.
\end{equation}
If neither coefficient on the right-hand side of \eqref{Eq product first obstruction} vanishes identically, choose $x_i\in M_i$ with $\widetilde a_{n_i+1}(g_i)(x_i)\neq0$, independently for $i=1,2$. Their product is nonzero at $(x_1,x_2)$. Conversely, if either coefficient vanishes identically, then \eqref{Eq BLF product condition} holds. This proves the equivalent assertion.

Finally, if both $\widetilde a_{n_i+1}(g_i)$ vanish identically, then $\widetilde a_j(g_i)\equiv0$ for every $j>n_i$. For $j>n$, each summand in the second formula of \eqref{Eq symmetric products} has either $\ell>n_1$ or $j-\ell>n_2$, and hence vanishes. This proves the last conclusion.
\end{proof}

\begin{rmk}
The distinction between
Theorems~\ref{Thm SLF product flatness} and~\ref{Thm BLF product flatness} can be seen intuitively as follows. At a fixed point, the formal twisted expansion of a Bergman logarithmically flat metric of dimension $n_i$ may contain a term $\widetilde a_{n_i+1}(g_i)m^{-1}$. Multiplying two such terms gives an $m^{-2}$ term. On a base of dimension $n=n_1+n_2$, this is exactly the coefficient with index $n+2$, which is incompatible with Bergman logarithmic flatness. This observation alone would only identify a possible obstruction; Remark~\ref{Rmk BLF product counterexample} shows that it actually occurs.
\end{rmk}

\subsection{Hermitian symmetric metrics}\label{Subsec Hermitian symmetric}
We prove Corollary~\ref{Cor Hermitian symmetric flatness} in this subsection. As recalled in \S~\ref{Sec Introduction}, Engli\v{s}--Zhang \cite{EnZh} established the vanishing of logarithmic terms in the Szeg\H{o} kernels of circle bundles and the Bergman kernels of disk bundles associated with compact Hermitian symmetric spaces. Related results over bounded symmetric domains, bounded homogeneous domains, and Euclidean spaces were obtained by Roos \cite{Roos05}, Ahn--Park \cite{AP12}, Hao--Wang \cite{HW15}, Ishi--Park--Yamamori \cite{IPY17}, and Yamamori \cite{Yam13}; see also the references therein. We combine these earlier results with Theorems~\ref{Thm SLF product flatness} and~\ref{Thm BLF product flatness}  to prove the corollary for arbitrary locally Hermitian symmetric K\"ahler metrics. 
	
\begin{proof}[Proof of Corollary~\ref{Cor Hermitian symmetric flatness}]
	We first consider the irreducible symmetric factors. By \eqref{Eq symmetric scaling}, we may choose a convenient positive normalization of each invariant metric.
	
	Let $(X,g)$ be an irreducible Hermitian symmetric space of compact type, of complex dimension $d$. Choose the homogeneous positive line bundle $L$ introduced at the beginning of \cite[Section~3]{EnZh}, and normalize $g$ so that $L$ polarizes it. Engli\v{s}--Zhang \cite[Lemma~3.1]{EnZh} give
	\[
		\mathcal B_m=\frac{P(m)}{\operatorname{Vol}_g(X)}.
	\]
	Here $P(m)$ is the scalar prefactor in \cite[Lemma~3.1]{EnZh}; the factor $\operatorname{Vol}_g(X)^{-1}$ accounts for replacing their invariant probability measure by $dV_g$. Moreover, $P$ is a polynomial of degree $d$, as stated explicitly in \cite[proof of Theorem~3.2]{EnZh}.
	
	At the beginning of \cite[Section~3]{EnZh}, the authors also identify $L^p$ with $\mathcal C_X^{-1}$, where $p$ is the genus. The invariant Hermitian metrics agree up to a positive constant, which does not affect the Bergman kernel function. Consequently,
		\[
		\widetilde{\mathcal B}_m=\mathcal B_{m-p}
		=\frac{P(m-p)}{\operatorname{Vol}_g(X)},\qquad m\geq p.
		\]
	Thus both kernel functions are polynomials of degree $d$. Then  \eqref{Eq two Bergman expansions} and uniqueness of the asymptotic coefficients yield $a_j(g)=\widetilde a_j(g)=0$ for $j>d$.
		
Next let $(X,g)$ be an irreducible Hermitian symmetric space of noncompact type, realized as a bounded symmetric domain $\Omega\subset\mathbb C^d$ via the Harish-Chandra realization.  Let $N(z,w)$ and $\gamma$ denote the generic norm and genus of $\Omega$, with $N(0,0)=1$. By \eqref{Eq symmetric scaling}, we may normalize $g$ so that it has K\"ahler potential $\varphi(z)=-\log N(z,z)$. Write $L=\Omega\times\mathbb C$ for the trivial holomorphic line bundle, with constant holomorphic frame $e$. Equip $L$ with the Hermitian metric
		\[
		h(e,e)(z)=e^{-\varphi(z)}=N(z,z).
		\]
Then $(L,h)$ polarizes $g$. Write $dV_0$ for the Euclidean volume form, and $V$ for the volume of $\Omega$ with respect to $dV_0$. The Bergman kernel of $\Omega$ is $V^{-1}N(z,z)^{-\gamma}$ (cf. \cite{Mok89, Roos05} and references therein). This kernel and $\det g$ obey the same transformation law under automorphisms of $\Omega$. Their quotient is therefore invariant and, by transitivity, constant (cf. \cite{Mok89}). Consequently, for some $C>0$,
		\[
		\det g=C N(z,z)^{-\gamma},\qquad dV_g=C N(z,z)^{-\gamma}\,dV_0.
		\]
	For $\beta>-1$, consider the Hilbert space of holomorphic functions
		\[
		\mathcal H_\beta=\left\{f\in\mathcal O(\Omega):\|f\|_\beta^2:=\int_\Omega |f(z)|^2N(z,z)^\beta\,dV_0(z)<\infty\right\}.
		\]
	Let $K_\beta(z,w)$ denote its reproducing kernel. Explicitly, if $\{u_{\beta,\nu}\}_\nu$ is an orthonormal basis of $\mathcal H_\beta$, then
		\[
		K_\beta(z,w)=\sum_\nu u_{\beta,\nu}(z)\overline{u_{\beta,\nu}(w)}.
		\]
	Roos \cite[Corollary~2.2, equation~(2.6)]{Roos05} gives
\begin{equation}\label{Eq symmetric scalar weighted kernel}
	K_\beta(z,z)=\frac{\chi(\beta)}{\chi(0)V}N(z,z)^{-\beta-\gamma},\qquad \beta>-1,
\end{equation}
	where the Hua polynomial $\chi$ has degree $d$ by \cite[Theorem~2.1, equation~(2.4)]{Roos05} (see also Ahn--Park \cite[Corollary~2.2]{AP12}). 

We now identify the two section spaces defining $\mathcal B_m$ and $\widetilde{\mathcal B}_m$ with the above weighted spaces of functions. Write $dz=dz_1\wedge\cdots\wedge dz_d$ for the standard holomorphic frame of the canonical bundle $\mathcal C_\Omega$. Its induced Hermitian metric satisfies $|dz|_{\mathcal C_\Omega}^2=(\det g)^{-1}$. 
Write holomorphic sections of $L^m$ and $L^m\otimes\mathcal C_\Omega$ as $f e^{\otimes m}$ and $q e^{\otimes m}\otimes dz$, respectively, where $f,q\in\mathcal O(\Omega)$. We have
\[
\|f e^{\otimes m}\|^2=\int_\Omega |f|^2N(z,z)^m\,dV_g=C\int_\Omega |f|^2N(z,z)^{m-\gamma}\,dV_0.
\]
\[
\|q e^{\otimes m}\otimes dz\|^2=\int_\Omega |q|^2N(z,z)^m(\det g)^{-1}\,dV_g=\int_\Omega |q|^2N(z,z)^m\,dV_0.
\]
For positive integers $m>\gamma-1$, we may therefore take
\[
\left\{C^{-1/2}u_{m-\gamma,\nu}e^{\otimes m}\right\}_\nu,\qquad \left\{u_{m,\nu}e^{\otimes m}\otimes dz\right\}_\nu
\]
as orthonormal bases of the Bergman spaces of $L^m$ and $L^m\otimes\mathcal C_\Omega$, respectively.
Summing their pointwise squared norms and applying \eqref{Eq symmetric scalar weighted kernel}, we obtain
	\[
	\mathcal B_m(z)=C^{-1}N(z,z)^m\sum_\nu|u_{m-\gamma,\nu}(z)|^2=C^{-1}N(z,z)^m K_{m-\gamma}(z,z)=\frac{\chi(m-\gamma)}{C\chi(0)V}.
	\]
	\[
	\widetilde{\mathcal B}_m(z)=N(z,z)^m(\det g(z))^{-1}\sum_\nu|u_{m,\nu}(z)|^2=C^{-1}N(z,z)^{m+\gamma}K_m(z,z)=\frac{\chi(m)}{C\chi(0)V}.
	\]
	The metric $g$ is complete and homogeneous, so its Ricci curvature is bounded below relative to $g$. The global kernel functions therefore have the expansions with coefficients $a_j(g)$ and $\widetilde a_j(g)$ by Proposition~\ref{Prop TYZ expansion} and \cite[Proposition~1.5]{EXX-BLF}, respectively. Since both kernel functions are polynomials in $m$ of degree $d$, uniqueness of the asymptotic coefficients gives $a_j(g)=\widetilde a_j(g)=0$ for every $j>d$.
	
	For Euclidean space $(\mathbb C^d,g_0)$, let $L=\mathbb C^d\times\mathbb C$ be the trivial holomorphic line bundle, with constant holomorphic frame $e$, and equip it with the Hermitian metric $h(e,e)(z)=e^{-|z|^2}$. Then $(L,h)$ polarizes $g_0$. Write $dz=dz_1\wedge\cdots\wedge dz_d$ for the standard holomorphic frame of the canonical bundle. Since $\det g_0=1$, we have $|dz|_{\mathcal C_{\mathbb C^d}}^2=1$ and $dV_{g_0}=dV_0$, where $dV_0$ is Euclidean volume. Consequently, for every positive integer $m$, the untwisted and canonically twisted section norms satisfy
	\[
	\|f e^{\otimes m}\|^2=\|f e^{\otimes m}\otimes dz\|^2 =\int_{\mathbb C^d}|f(z)|^2e^{-m|z|^2}\,dV_0(z).
	\]
Thus both Bergman spaces of sections can be identified with the Fock space, where the latter consists of entire functions square integrable with respect to $e^{-m|z|^2}\,dV_0$.  By the classical formula for the reproducing kernel $K_m(z,w)$ of this Fock space (cf.~\cite{Yam13}),
\[
K_m(z,w)=\frac{m^d}{\pi^d}e^{m\langle z,w\rangle},\qquad \text{where}~~ \langle z,w\rangle=\sum_{j=1}^d z_j\overline{w_j}.
\]
Choose an orthonormal basis $\{f_\nu\}_\nu$ of the Fock space. Then $\{f_\nu e^{\otimes m}\}_\nu$ and $\{f_\nu e^{\otimes m}\otimes dz\}_\nu$ are orthonormal bases of the two Bergman spaces, respectively. Hence
\[
\mathcal B_m(z)=\widetilde{\mathcal B}_m(z)=e^{-m|z|^2}\sum_\nu|f_\nu(z)|^2=e^{-m|z|^2}K_m(z,z)=\frac{m^d}{\pi^d}.
\] Comparing this with \eqref{Eq two Bergman expansions}, we obtain $a_j(g_0)=\widetilde a_j(g_0)=0$ for every $j\geq1$.

Finally, a locally Hermitian symmetric K\"ahler metric is locally holomorphically isometric to a product of a flat factor and irreducible Hermitian symmetric factors of compact or noncompact type; see \cite{Mok89,Hel01}. Write the factors as $(M_i,g_i)$, with $n_i=\dim_{\mathbb C}M_i$ and $n=\sum_i n_i$. The preceding arguments and \eqref{Eq symmetric scaling} give $a_j(g_i)=\widetilde a_j(g_i)=0$ for $j>n_i$. 
Repeated application of \eqref{Eq symmetric products} yields
	\[
	a_j(g_1\oplus\cdots\oplus g_r)=\sum_{j_1+\cdots+j_r=j}\prod_{i=1}^r a_{j_i}(g_i),\qquad \widetilde a_j(g_1\oplus\cdots\oplus g_r)=\sum_{j_1+\cdots+j_r=j}\prod_{i=1}^r\widetilde a_{
	j_i}(g_i),
	\]
where all indices are nonnegative. If $j>n$, every summand has $j_i>n_i$ for at least one $i$, and hence vanishes. 
Since the coefficients are invariant under local holomorphic isometries, \eqref{Eq symmetric coefficient vanishing} holds on $M$. The assertions about Szeg\H{o} and Bergman logarithmic flatness follow from \eqref{Eq comparison criteria}.
\end{proof}

\def\cprime{$'$}

\end{document}